\documentclass{article}

\date{}

\usepackage{amssymb}
\usepackage{epsfig}

\newcommand\eq[1] {(\ref{#1})}

\newcommand{\bfm}[1]{\mbox{\boldmath ${#1}$}}

\newcommand{\beqa}{\begin{eqnarray}}
\newcommand{\eeqa}{\end{eqnarray}}
\newcommand{\bequ}{\begin{equation}}
\newcommand{\eequ}[1]{\label{#1}\end{equation}}
\newcommand{\Grad}{\nabla}

\newcommand{\Ga}{\alpha}

\newcommand{\Gd}{\delta}

\newcommand{\Gve}{\varepsilon}

\newcommand{\Gs}{\sigma}

\newcommand{\Go}{\omega}

\newcommand{\GD}{\Delta}

\newcommand{\GO}{\Omega}

\newcommand{\BGx}{\bfm\xi}

\newcommand{\CK}{{\cal K}}
\newcommand{\CL}{{\cal L}}
\newcommand{\CM}{{\cal M}}
\newcommand{\CN}{{\cal N}}

\newcommand{\CQ}{{\cal Q}}
\newcommand{\CR}{{\cal R}}

\newcommand{\BCC}{{\bfm{\cal C}}}

\def\Bx{{\bf x}}
\def\By{{\bf y}}
\def\Bz{{\bf z}}

\def\BC{{\bf C}}
\def\BD{{\bf D}}

\def\BI{{\bf I}}

\def\BO{{\bf O}}

\def\BS{{\bf S}}

\def\BV{{\bf V}}

\def\BX{{\bf X}}
\def\BY{{\bf Y}}

\newcommand{\beq}{\begin{equation}}
\newcommand{\eeq}{\end{equation}}
\newcommand{\overliner}{\begin{eqnarray}}
\newcommand{\earr}{\end{eqnarray}}
\newcommand{\beqn}{\begin{equation*}}
\newcommand{\eeqn}{\end{equation*}}
\newcommand{\overlinern}{\begin{eqnarray*}}
\newcommand{\earrn}{\end{eqnarray*}}
\newcommand{\prt}{\partial}

\newcommand{\fr}{\frac}

\def\l{\label}

\begin{document}

\title{\Large \bf Asymptotic treatment of perforated domains  without homogenization }

\author
{\bf V. Maz'ya$^1$, 
A. Movchan$^2$  
\\ \\
{\small $^1$ Department of Mathematical Sciences, University of Liverpool,
 Liverpool L69 3BX,
U.K., } \\
{\small and  Department of Mathematics, Link\"oping University,  SE-581 83 Link\"oping, Sweden }\\
{\small $^2$ Department of Mathematical Sciences, University of Liverpool, Liverpool L69 3BX, U.K. }}


\maketitle

\vspace{.2in}

\centerline{\it To the memory of  Erhard Schmidt}

\vspace{.2in}

\begin{abstract}
{As a main result of the paper, we construct and justify an asymptotic approximation of Green's function in a domain with many small inclusions. Periodicity of the array of inclusions 
is not required.
We start with an analysis of the Dirichlet problem for the Laplacian in such a domain to  illustrate a 
method of meso scale asymptotic approximations for solutions of boundary value problems in multiply perforated domains. 
The asymptotic formula obtained involves a linear combination of solutions to certain model problems whose coefficients satisfy a linear algebraic system.
The solvability of this system is proved under weak geometrical assumptions, and  both uniform and energy estimates for the remainder term are derived. 

In the second part of the paper, the method is applied to derive an asymptotic representation 
of the Green's function in the same perforated domain. 
The important feature is the uniformity of the remainder estimate with respect to the independent variables.}
\end{abstract}

\noindent {\bf Keywords:} {\small Singular perturbations, 
meso scale approximations, multiply perforated domains, Green's function }

\section{Introduction}
\label{intro}

Uniform
asymptotic approximations of Green's kernels for various singularly and
regularly perturbed domains were constructed in \cite{CRAS}--\cite{MMN}. In particular, the papers \cite{CRAS},  \cite{JCAM} 
address the case of domains
containing several small inclusions with different types of boundary conditions.
In the present paper, 
a similar geometrical configuration is considered, but the number of inclusions becomes  a large parameter, which makes asymptotic formulae in  \cite{CRAS}, \cite{JCAM} 
inapplicable. 

In Sections \ref{formal_as}--\ref{energy_est}, we address the Dirichlet problem for the Poisson equation $-\GD u = f$ in a multiply perforated domain with zero Dirichlet data on the boundary. Section
\ref{formal_as} contains the formal asymptotic representation
\beq
u(\Bx) \sim v_f(\Bx) +   \sum_{j=1}^N C_j \Big(P^{(j)}(\Bx)   - 4 \pi ~\mbox{cap} (F^{(j)}) ~ H(\Bx, \BO^{(j)})                 \Big) ,
\eequ{intro_1}
where
\begin{itemize}
\item  $v_f$ is the solution of the same equation in a domain $\GO$ without inclusions,
\item $P^{(j)}$ is the harmonic capacitary potential of the inclusion $F^{(j)}$,
\item $\mbox{cap}(F^{(j)})$ is the  harmonic capacity of $F^{(j)}$,
\item $H$ is the regular part of Green's function $G$ 
of $\GO$. 
\end{itemize}
The coefficients $C_j$ should be found from the algebraic system
\beq
 (\BI + \BS \BD) \BC + \BV_f = {\bf 0},
\eequ{intro_2}
where $\BI$ is the identity matrix, and the matrices $\BS, \BD$ and the vectors $\BC, \BV_f$ are defined by
 \beq
 \BS  
 = \Big\{
 (1-\Gd_{ik}) G(\BO^{(k)}, \BO^{(i)})
 \Big\}_{i,k=1}^N,
 ~
\BD =  4 \pi ~\mbox{diag} ~ \{ \mbox{cap} (F^{(1)}), \ldots,
\mbox{cap} (F^{(N)}) \},
\eequ{intro_3}
and
\beq
\BC=(C_1, \ldots, C_N)^T, ~
 \BV_f = (v_f(\BO^{(1)}), \ldots, v_f(\BO^{(N)}))^T,
 \eequ{intro4}
 with $\BO^{(j)}$ being interior points in $F^{(j)}.$
 
 The unique solvability of the system \eq{intro_2} is not obvious, and it is established in Section \ref{aux_syst} under the natural assumption $\Gve < c~ d$, where $c$ is a sufficiently small absolute constant, $\Gve$ is the maximum of diameters of $F^{(j)}, j = 1, \ldots, N$, and $d$ characterizes the distance between inclusions.  In the same section, we obtain auxiliary estimates for the vector $\BC$ under two different constraints on $\Gve$ and $d$. 
 
 Justification of the formal asymptotic approximation \eq{intro_1} is given in Sections \ref{uniform} and \ref{energy_est}, where we show that the remainder admits the uniform estimate $O(\Gve + \Gve^2 d^{-7/2})$ and
 the energy estimate $O(\Gve^2 d^{-4})$. 
 
 Although we see that the asymptotic method described above may be applied formally under a very mild geometrical constraint $\Gve < c ~ d$, the convergence of the approximation in the space  $L_\infty(\GO)$ and
 in the Sobolev space $H^1(\GO)$ has been proved when $\Gve \ll d^{7/4}$ and $\Gve \ll d^2,$ respectively. Hence, the asymptotic approximation \eq{intro_1} is efficient for certain meso scale geometries, intermediate between a collection of inclusions whose size $\Gve$ is comparable with $d$ and the classical situation with $\Gve \sim \mbox{const}~ d^3$ appearing in the homogenization theory (see \cite{Mar}, \cite{Cior1} et al.).
 As is well known, in the latter case $u$ is approximated by a solution of the equation with a ``strange term'' $-\GD \hat{u}+ \mu \hat{u} =f, ~ \mu \geq 0.$
 
 In the concluding Section \ref{meso_green}, we derive the above mentioned asymptotic formula for Green's function $G_N(\Bx, \By)$, uniform with respect to $\Bx$ and $\By$. The following is a specially simple form in the case of $\GO = {\Bbb R}^3:$
 $$
 G_N(\Bx, \By) = \fr{1-N}{4 \pi |\Bx - \By|} + \sum_{j=1}^N g^{(j)}(\Bx, \By) + \sum_{1 \leq i,j \leq N, ~ i \neq j} {\cal C}_{ij} P^{(i)}(\Bx) P^{(j)}(\By) + O(\Gve d^{-2}),
 $$
 where ~~$g^{(j)}$ ~~are ~~ Green's ~~functions ~~in ~~${\Bbb R}^3 \setminus F^{(j)}$, ~~and~~ the~~ matrix~~ $\BCC =({\cal C}_{ij})_{i,j=1}^N$ ~~is 
 defined by $\BCC = (\BI + \BS \BD)^{-1} \BS.$

\section{Main notations and formulation of the problem in the perforated region} 
\l{notations}

Let $\GO$ be an arbitrary domain in ${\Bbb R}^3,$  
and let $\{  \BO^{(j)}\}_{j=1}^N
$ and $\{F^{(j)}\}_{j=1}^N$ be collections of points and disjoint compact subsets of $\GO$ such that
$\BO^{(j)} \in F^{(j)},$ and $F^{(j)}$ have positive harmonic capacity. 
Assume that the diameter $\Gve_j$ of $F^{(j)}$ is small compared to the diameter of 
$\GO$. 
 We 
 shall also use the notations
\beq
d= 
2^{-1} \min_{i \neq j, 1 \leq i, j \leq N} |\BO^{(j)} - \BO^{(i)}|, ~~ \Gve = \max_{1 \leq j \leq N} ~ \Gve_j. 
\eequ{eps_d}
It is assumed that $\Gve < c ~d,$ with $c$ being a sufficiently small constant.

We require that there exists an open set $\Go$ 
such that
\beq
\bigcup_{j=1}^N F^{(j)} \subset \Go, ~ \mbox{diam}(\Go) = 1, ~
 \mbox{dist}~ (\prt \Go, \prt \GO)  \geq 
 ~  2 d, 
~~~~\mbox{and      }~~~~ \mbox{dist}  \Big\{  \bigcup_{j=1}^N F^{(j)}, \prt \Go  \Big\} \geq 
 ~ 2 d.
\eequ{dist}

\begin{figure}[h]
\begin{center}
\includegraphics[scale=0.3]{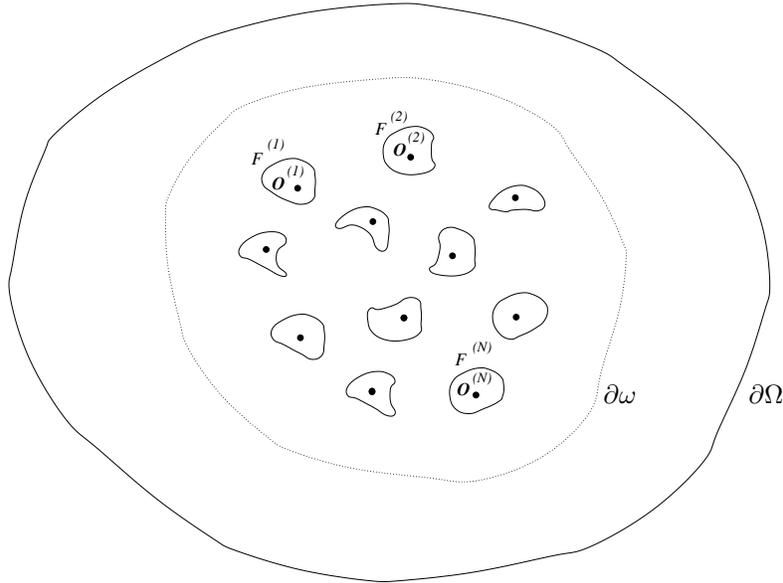} 
\begin{picture}(0,0)(0,0)

\put(-15,68){ $\prt \GO$}

\put(-70,68){ $\prt \Go$}

\end{picture}
\end{center}
\caption{Perforated domain containing many holes.} \label{fig1}
\end{figure}

Let us introduce the complimentary domain
 \beq
 \GO_{N} = \GO \setminus \cup_{j=1}^N
 F^{(j)},
 \eequ{domain}
 as shown in Fig. \ref{fig1}.

 Let $u$ denote the variational solution of the Dirichlet problem
 \beqa
 -\GD u
 (\Bx) &&=  f(\Bx) , ~~ \Bx \in \GO_{N}, \l{f1}
 \\
 u
 (\Bx) &&= 0, ~~ \Bx \in \prt \GO_{ N}, \l{f2}
 \end{eqnarray}
 where $f$ is assumed to be a smooth function 
  with a compact support in $\GO$, such that
  $\mbox{diam}(\mbox{supp}~ f) \leq C$ with $C$ being an absolute constant.  

We seek an asymptotic approximation of $u
$ as $N \to \infty.$


\section{Auxiliary problems}
\l{aux_pr}

We collect here solutions of some boundary value problems to be used in the asymptotic approximation of $u$. 

\subsection{Solution of the unperturbed problem}
 \l{unpert_problem}

By $v_f$ we mean the  variational solution of the Dirichlet problem
 \beqa
 -\GD v_{f}(\Bx)  &&= f(\Bx), ~~ \Bx \in \GO, 
 \l{f1u}
 \\
 v_{f}(\Bx) &&= 0, ~~ \Bx \in \prt \GO, 
 \l{f2u}
 \end{eqnarray}
 where $f$ is the same smooth function as in \eq{f1}.

 \subsection{Capacitary potentials of $F^{(j)}$}
 \l{model_fields}

The harmonic capacitary potential of $F^{(j)}$ will be
denoted by $P^{(j)}$, and it is defined
as a unique variational solution of the Dirichlet problem
 \beqa
 \GD P^{(j)}(\Bx) &&= 0 ~~\mbox{on} ~~ {\Bbb R}^3 \setminus F^{(j)},
 \l{eq6}\\
 P^{(j)} (\Bx) &&= 1~~\mbox{for} ~~ \Bx \in \prt ({\Bbb R}^3 \setminus
 F^{(j)}), \l{eq7} \\
 P^{(j)} (\Bx) &&= O(\Gve
 |\Bx-\BO^{(j)}|^{-1})
  ~~ \mbox{as} ~~ \Gve
  ^{-1} |\Bx-\BO^{(j)}| \to \infty.
 \l{eq8}
 \eeqa
It is well known (see, for example, \cite{PS}),
that  these functions have the following asymptotic
representations:
  \beq
  P^{(j)}(\Bx) = \fr{\mbox{\rm cap}(F^{(j)})}{
 |\Bx - \BO^{(j)}|} + O(\Gve
~ \mbox{cap} (F^{(j)}) |\Bx-\BO^{(j)}|^{-2}) ~~
 \mbox{for}~~ |\Bx - \BO^{(j)}| > 2 \Gve . 
 \eequ{P}
The harmonic capacity of
the set $F^{(j)}$ can be found by
 \beq
 \mbox{\rm cap}(F^{(j)}) = \fr{1}{4 \pi}\int_{{\Bbb R}^3 \setminus F^{(j)}}
 |\Grad P^{(j)}(\BGx)|^2 d
 \BGx.
 \eequ{cap}

 \subsection{Green's function for the unperturbed domain}
 \label{Green_unpert}

Green's function for the unperturbed domain  is denoted
 by $G(\Bx, \By)$, and it satisfies the boundary value problem
  \beqa
  {\GD}_x G(\Bx, \By) + \Gd(\Bx - \By) = 0, ~~ \Bx, \By \in \GO, \l{eq9}
  \\
  G(\Bx, \By) = 0 ~~\mbox{as} ~~\Bx \in \prt \GO ~~ \mbox{and} ~~\By
  \in \GO. \l{eq10}
  \eeqa
  The  regular part of Green's function is defined by
  \beq
  H(\Bx, \By) = 
  ({4 \pi |\Bx - \By|})^{-1} - G(\Bx, \By).
  \eequ{H}

\section{Formal asymptotic algorithm}
\l{formal_as}

Let the solution $u$ of \eq{f1},  \eq{f2} be written as
\beq
u(\Bx) = v_f(\Bx) + R^{(1)}(\Bx),
\eequ{falg_1}
where $v_f$ 
solves the auxiliary Dirichlet problem \eq{f1u}, \eq{f2u} in the unperturbed domain, whereas the 
function $R^{(1)}$ is harmonic in $\GO_N$ and satisfies the boundary conditions
\beq
R^{(1)}(\Bx) = 0~~\mbox{when} ~~ \Bx \in \prt \GO,
\eequ{falg2}
and
\beq
R^{(1)}(\Bx) = - v_f(\Bx) = - v_f(\BO^{(k)}) + O(\Gve) ~~ \mbox{when} ~~ \Bx \in \prt ( {\Bbb R}^3 \setminus F^{(k)} ) .
\eequ{falg3}
 Let us approximate the function $R^{(1)}$ in the form
 \beq
 R^{(1)}(\Bx) \sim \sum_{j=1}^N C_j \Big(P^{(j)}(\Bx)   - 4 \pi ~\mbox{cap} (F^{(j)}) ~ H(\Bx, \BO^{(j)})                 \Big),
\eequ{falg4}
where $C_j$ are unknown constant coefficients, and $P^{(j)}$ and $H$ are the same as in \eq{eq6}--\eq{P} and \eq{H}, respectively.

By \eq{P}, \eq{H} and \eq{eq10}, we deduce
\beq
P^{(j)}(\Bx)   - 4 \pi ~\mbox{cap} (F^{(j)}) ~ H(\Bx, \BO^{(j)})  = O(\Gve ~\mbox{cap} (F^{(j)})  |\Bx - \BO^{(j)}|^{-2}  ),
\eequ{falg5}
for all $\Bx \in \prt \GO, ~ j = 1, \ldots, N.$ 

On the boundary of a small inclusion $F^{(k)}~~$ ($k=1,\ldots,N$) 
we have
\beq
v_f(\BO^{(k)}) + O(\Gve) + C_k (1 + O(\Gve)) 
\eequ{falg6}
$$
+ \sum_{1\leq j \leq N, ~ j \neq k}C_j \Big(   4 \pi ~ \mbox{cap} (F^{(j)})~ G(\BO^{(k)}, \BO^{(j)})  + O(\Gve ~\mbox{cap} (F^{(j)}) | \Bx - \BO^{(j)}|^{-2}  )        \Big) = 0,
$$
for all $\Bx \in \prt ( {\Bbb R}^3 \setminus F^{(k)} ).$

Equation \eq{falg6} suggests that 
the constant coefficients $C_j, ~ j=1,\ldots, N,$ 
should be chosen to satisfy the system of linear algebraic equations
\beq
v_f(\BO^{(k)}) + C_k 
+ 4 \pi \sum_{1\leq j \leq N, ~ j \neq k}C_j  ~ \mbox{cap} (F^{(j)})~ G(\BO^{(k)}, \BO^{(j)})         = 0,
\eequ{falg7}
where $k =1,\ldots, N.$

Then within certain constraints on the small parameters $\Gve$ and $d$  (see \eq{eps_d}), it will be shown in the sequel that the above system of algebraic equations is solvable and that
the harmonic function
$$
R^{(2)}(\Bx) =  R^{(1)}(\Bx) - \sum_{j=1}^N C_j \Big(P^{(j)}(\Bx)   - 4 \pi ~\mbox{cap} (F^{(j)}) ~ H(\Bx, \BO^{(j)})                 \Big)
$$
is small on $\prt \GO_N$. Further application of the maximum principle for harmonic functions leads to an estimate of the remainder $R^{(2)}$ in $\GO_N.$

Hence, the solution \eq{falg_1} takes the form
\beq
u(\Bx) = v_f(\Bx) +   \sum_{j=1}^N C_j \Big(P^{(j)}(\Bx)   - 4 \pi ~\mbox{cap} (F^{(j)}) ~ H(\Bx, \BO^{(j)})                 \Big)  + R^{(2)}(\Bx),
\eequ{falg_8}
where $C_j$ are obtained from the algebraic system \eq{falg7}.

\section{
Algebraic system}
\l{aux_syst}

In this section we analyse the solvability of the system \eq{falg7}, and subject to certain constraints on $\Gve$ and $d$, derive estimates for the coefficients $C_j, ~ j=1,\ldots, N.$

The following matrices $\BS$ and $\BD$ will be used here: 
\beq
 \BS  
 = \Big\{
 (1-\Gd_{ik}) G(\BO^{(k)}, \BO^{(i)})
 \Big\}_{i,k=1}^N,
 \eequ{S_matrix1}
 and
 \beq
\BD =  4 \pi ~\mbox{diag} ~ \{ \mbox{cap} (F^{(1)}), \ldots,
\mbox{cap} (F^{(N)}) \}.
 \eequ{D_matrix1}
If the matrix
$\BI + 
\BS \BD$  is non-degenerate, 
then the components of the column vector $\BC=(C_1, \ldots, C_N)^T$ are defined by
\beq
\BC= - (\BI + 
 \BS \BD)^{-1} \BV_f , 
\eequ{the_system_soln}
where \beq \BV_f = (v_f(\BO^{(1)}), \ldots, v_f(\BO^{(N)}))^T. \eequ{BV} 

 Prior to the formulation of the result on the uniform asymptotic approximation of the solution to problem \eq{f1}-\eq{f2},
 we formulate and prove 
 auxiliary statements incorporating the  invertibility of the matrix $\BI + \BS \BD$ and   
estimates for 
components of the vector \eq{the_system_soln}. 

{\bf Lemma 1.} {\em If $
\max_{1 \leq j \leq N} \mbox{\rm cap}(F^{(j)})< 5 d/(24 \pi),$  then the matrix $\BI + \BS \BD$ is invertible and the column vector $\BC$ in {\rm \eq{the_system_soln}} satisfies the estimate 
\beq
\sum_{j=1}^N \mbox{\rm cap}(F^{(j)}) ~ C_j^2 \leq 
(1- \fr{24 \pi 
}{5 d} \max_{1 \leq j \leq N} \mbox{\rm cap}(F^{(j)}))^{-2} 
\sum_{j=1}^N \mbox{\rm cap}(F^{(j)})~(v_f(\BO^{(j)}))^2.
\eequ{eq18a}
}

{\em Proof:} 

According to \eq{the_system_soln}, we have $(\BI + \BS \BD) \BC = -\BV_f.$ Hence
\beq
\langle \BC, \BD \BC \rangle + \langle  \BS \BD \BC, \BD \BC \rangle = - \langle   \BV_f, \BD \BC \rangle.
\eequ{eq**}
Obviously, the right-hand side in \eq{eq**} does not exceed 
\beq
\langle \BC, \BD \BC \rangle^{1/2}  \langle \BV_f, \BD \BV_f \rangle^{1/2}.
\eequ{eq***}
Consider the second term in the left-hand side of \eq{eq**}. Using the mean value theorem for harmonic functions we deduce
\beqa
&& \langle   \BS \BD \BC, \BD \BC \rangle = (4 \pi )^2 \sum_{i \neq j, 1 \leq i,j \leq N}  G(\BO^{(i)}, \BO^{(j)}) \mbox{cap} (F^{(i)})  \mbox{cap} (F^{(j)}) ~ C_i C_j \nonumber \\
&&= (4 \pi)^2 \sum_{i \neq j, 1 \leq i,j \leq N}  \fr{\mbox{cap} (F^{(i)})  \mbox{cap} (F^{(j)}) ~ C_i C_j}{
|B^{(i)}| ~ |B^{(j)}|} \int_{B^{(i)}}  \int_{B^{(j)}}  G(\BX, \BY) d \BX d \BY, \nonumber
\end{eqnarray}
where $B^{(j)}= \{ \Bx: |\Bx - \BO^{(j)}| < d \}, ~ j = 1, \ldots, N, $ are non-overlapping balls of radius $d$ with the centers at $\BO^{(j)}$, and $|B^{(j)}| = 4 \pi d^3/3$ are the volumes of the balls. 
Also, the notation $B_d$ is  used here for the ball of radius $d$ with the center at the origin.

Let 
$\Xi(\Bx)$ be a piecewise function defined on $\GO$ as 
$$\Xi(\Bx) = \left\{ \begin{array}{cc}  
C_j \mbox{cap} ( F^{(j)}) & \mbox{in } ~~ B^{(j)}, j = 1, \ldots, N, \\
0 & \mbox{otherwise} .
\end{array} \right. $$

Then
\beqa
&& \langle   \BS \BD \BC, \BD \BC \rangle = \fr{9}{d^6} \bigg( \int_{\GO}  \int_{\GO} G(\BX, \BY) \Xi(\BX) \Xi(\BY) d \BX d \BY   \nonumber \\
&& 
-\sum_{j=1}^N (\mbox{cap} ( F^{(j)}))^2 C_j^2 \int_{B^{(j)}} 
 \int_{B^{(j)}} 
 G(\BX, \BY)  d \BX d \BY \bigg).
\l{l_bound_a}
\end{eqnarray}
The first term in the right-hand side of \eq{l_bound_a} is non-negative, which follows from the relation
\beqa
\int_{\GO}  \int_{\GO} G(\BX, \BY) \Xi(\BX) \Xi(\BY) d \BX d \BY 
 = 
 \int_\GO \Big|   \Grad_X \int_\GO G(\BX, \BY) \Xi(\BY) d \BY \Big|^2 d \BX \geq 0.
\l{non_negative}
\end{eqnarray}
The integral
$$
\int_{B^{(j)}} 
\int_{B^{(j)}} 
G(\BX, \BY)  d \BX d \BY 
$$
in the right-hand side of \eq{l_bound_a} allows for the estimate
\beqa
&& \int_{B^{(j)}}  \int_{B^{(j)}} G(\BX, \BY)  d \BX d \BY \leq \fr{1}{4 \pi} \int_{B_d}  \int_{B_d} \fr{d \BX d \BY}{|\BX - \BY|} \nonumber \\
&& = \fr{1}{4 \pi}  \int_{B_d} d \BX \Big\{  \int_{|\BY | < |\BX  |} \fr{d \BY}{| \BX - \BY  |}  + \int_{d > |\BY | >  |\BX  |} \fr{d \BY}{| \BX - \BY  |}    \Big\} \nonumber \\
&& = \fr{1}{4 \pi}   \int_{B_d} d \BX \Big\{  \int_0^{|\BX|}  d \rho  \int_{\{ \BY: |\BY| = \rho \}} \fr{d S_Y}{| \BX - \BY |} +  \int_{|\BX|}^{d}  d \rho \int_{ \{\BY: |\BY| = \rho  \}} \fr{d S_Y}{| \BX - \BY |}     \Big\}.
\l{G_ball_est}
\end{eqnarray}
Using the mean value theorem for harmonic functions we deduce
\beq
 \int_{\{ \BY: |\BY| = \rho \}} \fr{d S_Y}{| \BX - \BY |} = 4 \pi \rho^2 |\BX|^{-1}~~\mbox{when}~ |\BX| > \rho.
\eequ{ball_eq_a}
On the other hand,
\beq
 \int_{\{ \BY: |\BY| = \rho \}} \fr{d S_Y}{| \BX - \BY |} = 4 \pi \rho~~\mbox{when}~ |\BX| < \rho,
\eequ{ball_eq_b}
which follows from the relation
\beqa
&& \int_{ \{ \BY:  |\BY| = \rho  \}  } \fr{d S_Y}{\rho | \BX - \BY |} = - \int_{ \{ \BY:  |\BY| = \rho  \}  }  \fr{\prt }{ \prt | \BY  |}   \fr{1}{| \BX - \BY  |} d S_Y \nonumber \\
&& = - \int_{ \{ \BY:  |\BY| <  \rho  \}  }  \GD_Y \fr{1}{| \BX - \BY  |}  d \BY = 4 \pi~~~ \mbox{when} ~~ |\BX| < \rho. \nonumber
\end{eqnarray}
It follows from \eq{G_ball_est}, \eq{ball_eq_a} and \eq{ball_eq_b} that
\beq
\int_{B^{(j)}} 
\int_{B^{(j)}} 
G(\BX, \BY)  d \BX d \BY \leq  \fr{1}{2}\int_{B_d}  \Big(   d^2 - \fr{|\BX|^2}{3}\Big) d \BX= \fr{8 \pi d^5}{15}.
\eequ{ball_G_est}
Next, \eq{eq**}, \eq{eq***}, \eq{l_bound_a} and \eq{ball_G_est} lead to 
\beqa
&& \langle   \BS \BD \BC, \BD \BC \rangle \geq \fr{9}{d^6} \int_{\Go} 
\int_{\Go} 
G(\BX, \BY) \Xi(\BX) \Xi(\BY) d \BX d \BY 
- 
 \fr{9  \Ga}{d} \sum_{j=1}^N (\mbox{cap} (F^{(j)}))^2 C_j^2 ,
\l{l_bound}
\end{eqnarray}
where
$
\Ga = 
\fr{8 \pi}{15}.$ 
Then \eq{eq**} and  \eq{l_bound} 
imply
\beqa
&& 
\Big(1- \fr{9 \Ga}{d}  \max_{1 \leq j \leq N} \mbox{cap} (F^{(j)}) \Big) \sum_{j=1}^N C_j^2 \mbox{cap} (F^{(j)})  \nonumber \\
&& 
\leq  
\Bigg({\sum_{j=1}^N C_j^2 \mbox{cap} (F^{(j)})}\Bigg)^{1/2}  
\Bigg({\sum_{j=1}^N (v_f(\BO^{(j)}))^2 \mbox{cap} (F^{(j)})} \Bigg)^{1/2},  \nonumber
\end{eqnarray}
which yields
\beq
(1- \fr{24 \pi 
}{5 d} \max_{1 \leq j \leq N} \mbox{cap} (F^{(j)})) 
\Bigg({\sum_{j=1}^N C_j^2 \mbox{cap} (F^{(j)}) }  \Bigg)^{1/2} \leq  
\Bigg({\sum_{j=1}^N (v_f(\BO^{(j)}))^2 \mbox{cap} (F^{(j)})}\Bigg)^{1/2}.
\eequ{est_lemma}
Thus, if $
\max_{1 \leq j \leq N} \mbox{cap} (F^{(j)}) < \fr{5}{24 \pi} d,$ then the matrix $\BI + \BS \BD$ is invertible and the estimate \eq{eq18a} holds. The proof is complete. $\Box$

Replacement of the inequality $\Gve < c d$  by the stronger constraint $ \Gve < c d^2$ 
leads to the statement  

{\bf Lemma 2.} {\em Let the small parameters $\Gve$ and $d$, defined in {\rm \eq{eps_d}}, satisfy
\beq
\Gve < c d^2,
\eequ{d2_constr}
where $c$ is a sufficiently small absolute constant.
Then the components $C_j$  of vector $\BC$ in {\rm \eq{the_system_soln}} allow for the estimate
\beq
|C_k| \leq c \max_{1 \leq j \leq N} |v_f(\BO^{(j)})|.
\eequ{L2_C}
}

{\em Proof:}  Let us write the system \eq{falg7} as
\beq
C_k + 4 \pi \sum_{1 \leq j \leq N, ~ j \neq k} C_j \fr{\mbox{cap} (F^{(j)})}{|B^{(j)}_{d/4}|} \int_{B^{(j)}_{d/4}} G(\BO^{(k)}, \By)  d \By = - v_f(\BO^{(k)}), ~ k=1, \ldots, N,
\eequ{L2_eq1}
where $B^{(j)}_{d/4}$ is the ball of radius $d/4$ with the centre at $\BO^{(j)}$.
Also let $\Gs$ be a piece-wise constant function such that
\beq
\Gs(\Bx) = \left\{    
\begin{array}{cc}
C_j ~\mbox{cap} (F^{(j)}), & \Bx \in B^{(j)}_{d/4} , \\
0, & \Bx \in {\Bbb R}^3 \setminus \cup_{m=1}^N B^{(m)}_{d/4}.
\end{array}
\right.
\eequ{L2_Gs}
Multiplying \eq{L2_eq1} by $\mbox{cap} ( F^{(k)})$ and writing the equations obtained in terms of $\Gs$ we get
$$
\Gs(\BO^{(k)}) + \fr{192}{d^3} ~ \mbox{cap} (F^{(k)})  \int_{\bigcup_{1 \leq j \leq N, ~ j \neq k} B^{(j)}_{d/4}} \Gs(\By) G(\BO^{(k)}, \By) d \By = -v_f(\BO^{(k)}) ~\mbox{cap} (F^{(k)}),
$$ 
which is equivalent to
\beq
\Gs(\BO^{(k)}) + \fr{192 
}{
 d^3} ~\mbox{cap} (F^{(k)}) ~ \int_{\bigcup_{1 \leq j \leq N}   B^{(j)}_{d/4}}     
G(\By,\Bz) {\Gs(\By) d \By}
 = \mbox{cap} (F^{(k)})~ \Phi^{(k)}(\Bz), ~ k = 1, \ldots, N,
\eequ{L2_eq2}
where
\beqa
&& \Phi^{(k)}(\Bz )= - v_f(\BO^{(k)}) + \fr{192 
}{
d^3} 
\int_{B^{(k)}_{d/4}} 
G(\By, \Bz) {\Gs(\By) d \By}
\l{L2_eq2a} \\
&& + \fr{192}{d^3} 
\int_{\bigcup_{1 \leq j \leq N, ~ j \neq k} B^{(j)}_{d/4}} \Gs(\By) \Big((H(\BO^{(k)}, \By) - H(\Bz, \By)\Big) d \By \nonumber \\
&&  + \fr{48}{\pi d^3} 
 \int_{\bigcup_{1 \leq j \leq N, ~ j \neq k} B^{(j)}_{d/4}} \Gs(\By) \Big\{    
 \fr{1}{| \By - \Bz |} - \fr{1}{ |  \BO^{(k)}   -  \By |   }
 \Big\} d \By, ~~ \mbox{for all}~ \Bz \in B^{(k)}_{d/4}. \nonumber 
 \end{eqnarray}
 Next, we multiply \eq{L2_eq2} by
 $$
 \Bigg(     \int_{\bigcup_{1 \leq j \leq N} B^{(j)}_{d/4}} 
 G(\By, \Bz) {\Gs(\By)}   d \By 
      \Bigg)^{2M-1},
 $$
 where $M
 $ is a positive 
 integer number. Also, taking into account that $\Gs(\BO^{(k)}) = \Gs(\Bz) ~~\mbox{for all}~ \Bz \in B^{(k)}_{d/4}$ we write
 \beqa
&&  \Bigg(     \int_{\bigcup_{1 \leq j \leq N} B^{(j)}_{d/4}} 
 G(\By, \Bz) {\Gs(\By) d \By} 
  \Bigg)^{2M-1} \Gs(\Bz) 
  + \fr{192 
  }{
   d^3} ~\mbox{cap} ( F^{(k)}) ~ \Bigg( \int_{\bigcup_{1 \leq j \leq N}   B^{(j)}_{d/4}}      
   G(\By, \Bz) {\Gs(\By) d \By}
   \Bigg)^{2M}  \nonumber \\
&&  = \mbox{cap} (F^{(k)})~ \Phi^{(k)}(\Bz)  \Bigg(     \int_{\bigcup_{1 \leq j \leq N} B^{(j)}_{d/4}} 
 G(\By, \Bz) {\Gs(\By) d \By}
       \Bigg)^{2M-1}, ~ \Bz \in 
 B^{(k)}_{d/4}. \nonumber
\end{eqnarray}
Since $\Gs =0$ outside the balls $B^{(k)}_{d/4}$, it follows that the integration of the above equation over $B^{(k)}_{d/4}$  and summation with respect to $k=1,\ldots,N$ lead to
\beqa
&& \int_{\GO 
} \Bigg(     \int_{ \GO 
} 
G(\By, \Bz) {\Gs(\By) d \By} 
      \Bigg)^{2M-1} \Gs(\Bz) d \Bz \nonumber \\
&& +  \fr{ 192 
}{
d^3}  \sum_{k=1}^N~\mbox{cap} ( F^{(k)}) ~ \int_{B^{(k)}_{d/4}} \Bigg( \int_{ \GO}  
G(\By, \Bz){\Gs(\By) d \By}
\Bigg)^{2M}  d \Bz \nonumber \\
&& = \sum_{k=1}^N\mbox{cap} (F^{(k)})~ \int_{B^{(k)}_{d/4}} \Phi^{(k)}(\Bz)  \Bigg(     \int_{  \GO }   
G(\By, \Bz){\Gs(\By) d \By}  
  \Bigg)^{2M-1} d \Bz.
\l{L2_eq3}
\end{eqnarray}
The identity
$$
\int_{ \GO 
} \Bigg(     \int_{    \GO 
} 
G(\By, \Bz){\Gs(\By) d \By}
     \Bigg)^{2M-1} \Gs(\Bz) d \Bz
= \fr{2M-1}{
M^2} \int_{\GO 
} \Bigg|  \Grad_z \Bigg( \int_{ \GO 
}      
G(\By, \Bz){\Gs(\By) d \By} 
\Bigg)^{M}  \Bigg|^2 d \Bz
$$
shows that the first term in the left-hand side of \eq{L2_eq3} is non-negative.  
By H\"{o}lder's inequality, the right-hand side of \eq{L2_eq3} does not exceed
\beqa
&& \Bigg(  \sum_{k=1}^N\mbox{cap} ( F^{(k)})~ \int_{B^{(k)}_{d/4}} (\Phi^{(k)}(\Bz))^{2M} d\Bz       \Bigg)^{1/(2M)} \nonumber \\
&& \times \Bigg(      \sum_{k=1}^N\mbox{cap} (F^{(k)})~ \int_{B^{(k)}_{d/4}}       \Bigg( \int_{ \GO  
}      
G(\By, \Bz){\Gs(\By) d \By} 
\Bigg)^{2 M} d \Bz         \Bigg)^{(2M-1)/(2M)}, \nonumber \\
\end{eqnarray}
and hence \eq{L2_eq3} yields
\beqa
&& \fr{192 
}{
 d^3}  \Bigg( \sum_{k=1}^N~\mbox{cap} (F^{(k)}) ~ \int_{B^{(k)}_{d/4}} \Bigg( \int_{  \GO 
 } 
G(\By, \Bz) {\Gs(\By) d \By} 
\Bigg)^{2M}  d \Bz \Bigg)^{1/(2M)} \nonumber \\
 &&  \leq \Bigg(  \sum_{k=1}^N\mbox{cap} (F^{(k)})~ \int_{B^{(k)}_{d/4}} (\Phi^{(k)}(\Bz))^{2M} d\Bz       \Bigg)^{1/(2M)}.
\end{eqnarray}
After the limit passage as $M \to \infty$ we arrive at
$$
d^{-3} \sup_{\Bz \in \bigcup_{1 \leq k \leq N} B^{(k)}_{d/4}} \Bigg|   \int_{ \GO 
 } 
 G(\By, \Bz) {\Gs(\By)    d\By}
  \Bigg| \leq c \max_{1 \leq k \leq N} ~ \sup_{\Bz \in B^{(k)}_{d/4}} |\Phi^{(k)}
(\Bz)|,
$$
and by \eq{L2_eq2}  we deduce 
\beq
|\Gs(\BO^{(k)})|  \leq c  ~\mbox{cap} (F^{(k)}) \max_{1 \leq j \leq N} \sup_{\Bz \in B^{(j)}_{d/4}} |\Phi^{(j)}(\Bz)|.
\eequ{L2_eq5}
In turn, it follows from the definition \eq{L2_eq2a} of the functions $\Phi^{(k)}$ that
\beqa
&& \sup_{\Bz  \in B^{(k)}_{d/4}} |\Phi^{(k)}
(\Bz)| \leq |v_f(\BO^{(k)})| 
+ \fr{192 
}{
d^3}   \max_{1 \leq q \leq N} |\Gs(\BO^{(q)})| 
 \sup_{\Bz  \in B^{(k)}_{d/4}}  \int_{B^{(k)}_{d/4}} G(\By, \Bz) 
{ d \By} 
\nonumber \\
&& +  \fr{192}{d^3}  \max_{1 \leq q \leq N} |\Gs(\BO^{(q)})|   \sup_{\Bz  \in B^{(k)}_{d/4}}   \sum_{1 \leq j \leq N, ~ j \neq k}
\int_{B^{(j)}_{d/4}} |H(\BO^{(k)}, \By)- H(\Bz, \By)| d \By \nonumber \\
 && +   \fr{48}{\pi d^3} \max_{1 \leq q \leq N} |\Gs(\BO^{(q)})|  
 \sup_{\Bz  \in B^{(k)}_{d/4}}   \sum_{1 \leq j \leq N, ~ j \neq k} \int_{ B^{(j)}_{d/4}}   
 \fr{|  \Bz - \BO^{(k)}  |}{| \By - \Bz |  |  \BO^{(k)}   -  \By |   }
d \By,  
 \nonumber 
 \end{eqnarray}
 which, together with \eq{L2_eq5}, 
 yields
 $$
 |\Gs(\BO^{(k)})| \leq c ~\mbox{cap} ( F^{(k)})~ \Big\{   \max_{1 \leq j \leq N}|  v_f(\BO^{(j)})  |  
 + d^{-2} \max_{1 \leq j \leq N} |  \Gs(\BO^{(j)})  |
 \Big\}.
 $$
 If  $\max_{1 \leq k \leq N} ~\mbox{cap} (F^{(k)}) < c d^2$, with $c$ being a sufficiently small constant, then referring to  the definition \eq{L2_Gs} of the function $\Gs$ we deduce \eq{L2_C},
 which completes the proof. $\Box$

\section{
Meso scale uniform 
approximation of $u$}
\l{uniform}


We obtain the next theorem, which is one of the principal results of the paper, 
under an additional assumption on the smallness of the capacities of $F^{(j)}$. 


{\bf Theorem 1.} {\em Let the 
parameters $\Gve$ and $d$, introduced in {\rm \eq{eps_d}}, satisfy the inequality
\beq
~~~~~~~~~~~~~~\Gve < c ~d^{7/4},
\eequ{eps_d_est}
where $c$ is a sufficiently small absolute constant.

Then the matrix $\BI + \BS \BD$, defined according to {\rm \eq{S_matrix1}, \eq{D_matrix1}},  is invertible, and the solution $u
(\Bx)$ to the
boundary value problem {\rm \eq{f1}--\eq{f2}} is defined by the
asymptotic formula
%
 \beq
 u
 (\Bx) =
v_f(\Bx) + \sum_{j=1}^N C_j \Big( P^{(j)}(\Bx) - 4 \pi
 ~\mbox{\rm
cap}(F^{(j)}) ~H(\Bx,
 \BO^{(j)}) \Big)  + R
(\Bx),
 \eequ{as_form}
where the column vector $\BC = (C_1, \ldots, C_N)^T$ is given by {\rm \eq{the_system_soln}}
 and
the remainder $R
(\Bx)$ is a function harmonic 
in
 $\GO_{N}$, which satisfies the estimate
\begin{eqnarray}
|R
(\Bx)| && \leq C 
\Big\{ \Gve \| \Grad v_f \|_{L_\infty (\Go)} + \Gve^2 d^{-7/2} \| v_f \|_{L_\infty (\Go)} \Big\} . 
\l{est_Th1} 
\end{eqnarray}
}


{\it  Proof:}    

The harmonicity of $R
$ follows directly from \eq{as_form}. 

If $\Bx \in \prt \GO$, then 
\beqa
R
(\Bx) =&& - 4 \pi \sum_{j=1}^N C_j \mbox{cap} ( F^{(j)}) \Big( \fr{1}{4 \pi |\Bx - \BO^{(j)}|} - H(\Bx, \BO^{(j)})\Big) 
\nonumber \\
&& + \sum_{j=1}^N |C_j| O(\Gve
~\mbox{cap} (F^{(j)}) |\Bx - \BO^{(j)}|^{-2}). \nonumber 
\end{eqnarray}
Since $G(\Bx, \BO^{(j)}) =0$ on $\prt \GO$, and $P^{(j)}$ satisfies \eq{P} we deduce
\beq
R
(\Bx)
= \sum_{j=1}^N O(\Gve
~\mbox{cap} ( F^{(j)})  |C_j|  |\Bx - \BO^{(j)}|^{-2}),
\eequ{Th1_Prf1}
where $|\Bx - \BO^{(j)}| \geq C~ d,$ and $C$ is a sufficiently large constant.

If $\Bx \in \prt (
{\Bbb R}^3 \setminus F^{(k)})$ then
\beqa
 R
(\Bx) =&& - v_f(\BO^{(k)}) + O(\Gve
\| \Grad v_f \|_{L_\infty(\Go)} )
+  
4 \pi \sum_{j=1}^N C_j \mbox{cap} (F^{(j)}) \Big( H(\BO^{(k)}, \BO^{(j)})  + O(\Gve
) \Big)
\nonumber \\
&& - C_k - \sum_{1 \leq j \leq N, ~ j \neq k } C_j \Big\{ \fr{\mbox{cap} ( F^{(j)}) }{| \BO^{(k)}- \BO^{(j)}|} + O(\fr{\Gve
~ \mbox{cap} ( F^{(j)}) }{|\BO^{(k)}- \BO^{(j)} |^2}) \Big\} .
\l{Th1_Prf2}
\end{eqnarray}
Noting that \eq{the_system_soln} can be written as the algebraic system
\beqa
&& C_k + 4 \pi \sum_{j=1}^N C_j (1-\Gd_{jk}) ~\mbox{cap} (F^{(j)}) \Big(   \fr{1}{4 \pi |\BO^{(k)}- \BO^{(j)}  |} 
\nonumber \\
&&
- H(\BO^{(k)}, \BO^{(j)}) \Big)  + v_f(\BO^{(k)}) =0,
\l{Prf1_system}
\end{eqnarray}
which, along with \eq{Th1_Prf2} and the obvious inequality $\mbox{cap} (F^{(j)}) \leq \Gve$, implies
\beqa
&& R
(\Bx) = O(\Gve
\| \Grad v_f \|_{L_\infty(\Go)}) + 4 \pi C_k \mbox{cap} (F^{(k)}) H( \BO^{(k)}, \BO^{(k)})
\nonumber \\
&& + \sum_{j=1}^N O(\Gve
~ \mbox{cap} (F^{(j)}) |C_j|) + \sum_{ 1 \leq  j \leq N, ~ j \neq k  } O(\fr{\Gve 
~\mbox{cap} ( F^{(j)}) }{| \BO^{(k)}- \BO^{(j)}  |^2} |C_j|) .
\l{Prf1a_system}
\end{eqnarray}

It suffices to estimate the sums
$$
\sum_{ 1 \leq j \leq N, ~ j \neq k} \fr{\Gve
~\mbox{cap} (F^{(j)}) |C_j|}{| \BO^{(k)}- \BO^{(j)} |^2}
$$
and 
$$
\sum_{1 \leq j \leq N} \fr{\Gve
~\mbox{cap} (F^{(j)})  |C_j|}{| \Bx- \BO^{(j)} |^2}, ~~ \Bx \in \prt \GO.
$$

When $\Gve < c ~d^{7/4}$ we refer to Lemma 1, and using the inequality \eq{eq18a} we derive 
\beqa
&&
\sum_{j\neq k, 1 \leq j \leq N}  \fr{\Gve 
~\mbox{cap} (F^{(j)})  |C_j|}{| \BO^{(k)}- \BO^{(j)} |^2}    \leq
\Bigg(\sum_{j \neq k, 1 \leq j \leq N} \fr{\Gve
^2~ \mbox{cap} (F^{(j)}) }{| \BO^{(k)}- \BO^{(j)} |^4}
\Bigg)^{1/2}  
\Bigg( \sum_{1 \leq j \leq N} \mbox{cap} ( F^{(j)}) C_j^2 
\Bigg)^{1/2}
\nonumber \\
&& \leq \mbox{const}~ d^{-1/2} 
\Bigg( \sum_{1 \leq j \leq N} \mbox{cap} (F^{(j)}) (v_f(\BO^{(j)}))^2 
\Bigg)^{1/2} 
\Bigg(\max_{1 \leq j \leq N} \Gve
^2 d
^{-3} \mbox{cap} (F^{(j)}) 
\Bigg)^{1/2}
\nonumber \\
&& \leq \mbox{const}~ \fr{\Gve^2}{d^{7/2}} \|v_f\|_{L_\infty (\Go)}.
\l{est1}
\end{eqnarray}
Similarly, when $\Bx \in \prt \GO$ we deduce
\beq
\sum_{1 \leq j \leq N} \fr{\Gve
~ \mbox{cap} (F^{(j)}) |C_j|}{| \Bx- \BO^{(j)} |^2} \leq  \mbox{const}~ \fr{\Gve^2}{d^{7/2}} \|v_f\|_{L_\infty (\Go)}.
\eequ{est2}

Combining \eq{est1}, \eq{est2}, \eq{Th1_Prf1} and \eq{Prf1a_system} we complete the proof  by referring to the classical maximum principle for harmonic functions.
~~$\Box$

Under the stronger constraint \eq{d2_constr} on $\Gve$ and $d$, Lemma 2 and representations  \eq{Th1_Prf1}, \eq{Th1_Prf2} lead to the following

{\bf Theorem  2.} {\it
If the inequality {\rm \eq{eps_d_est}} is replaced by  {\rm \eq{d2_constr}}, then the remainder term from {\rm \eq{as_form}} satisfies the estimate
\begin{eqnarray}
|R
(\Bx)| && \leq C 
\Big\{ \Gve \| \Grad v_f \|_{L_\infty (\Go)} + \Gve^2 d^{-3} \| v_f \|_{L_\infty (\Go)} \Big\} . 
\l{est_Cor} 
\end{eqnarray}
}

\section{The energy estimate}
\l{energy_est}

Under 
the constraint \eq{d2_constr} on $\Gve$ and $d$,  which is stronger than \eq{eps_d_est}, we derive the energy estimate for the remainder $R$.
This result is important, since it allows for the generalization to general elliptic systems, and in particular to elasticity where the classical maximum principle cannot be applied.  

{\bf Theorem 3.} {\it Let the 
parameters $\Gve $ and $d$, introduced in {\rm \eq{eps_d}}, satisfy the inequality
\beq
~~~~~~~~~~~~~~~~~~~~~~~~~~\Gve < c ~ d^2,
\eequ{eps_d_est_2}
where $c$ is a sufficiently small absolute constant.
Then the remainder $R
$ in {\rm \eq{as_form}} 
satisfies the estimate
\beq
\| \Grad R
\|_{L_2(\GO_N)} \leq \mbox{\rm Const} ~ \fr{\Gve^2}{d^4} \|f\|_{L_\infty(\GO_N)}
\eequ{th2_est}
}

{\bf Proof.} 


For every $k=1,\ldots,N,$ we introduce the function
\beqa
\Psi_k(\Bx) =&& v_f (\Bx)- v_f(\BO^{(k)})   + \sum_{1 \leq j \leq N, j\neq k} C_j \Big(P^{(j)}(\Bx) - 
\fr{
\mbox{cap} (F^{(j)})}{|\BO^{(k)}-\BO^{(j)}|}
\Big) 
\l{Psi_1} \\
&&- 4 \pi \sum_{j=1}^N
 C_j ~\mbox{cap} (F^{(j)}) \Big(  H(\Bx, \BO^{(j)})   - H(\BO^{(k)}, \BO^{(j)})  \Big)
\nonumber \\
 && - 4 \pi C_k ~\mbox{cap} ( F^{(k)}) H(\BO^{(k)}, \BO^{(k)}),
 \nonumber 
 \end{eqnarray}
where 
the coefficients $C_j$ satisfy the system \eq{falg7}.  

By \eq{as_form} 
and \eq{Psi_1}, for quasi-every $\Bx \in \prt ( {\Bbb R}^3 \setminus F^{(k)})$
\beqa
&&R(\Bx) + \Psi_k(\Bx) = 
- v_f(\BO^{(k)})  - C_k
\nonumber \\
&& -  \sum_{1 \leq j \leq N, ~ j \neq k} C_j \Big( 
\fr{
\mbox{cap} (F^{(j)}) }{| \BO^{(k)} - \BO^{(j)}  |}- 4 \pi  ~\mbox{cap} (F^{(j)}) H(\BO^{(k)}, \BO^{(j)}) \Big), \nonumber
\end{eqnarray}
which together with \eq{falg7} 
implies
$$
R(\Bx) + \Psi_k(\Bx) = 0
$$
quasi-everywhere on $\prt ({\Bbb R}^3 \setminus F^{(k)})$ (i.e. outside of a set with zero capacity).


The function $\Psi_0$, defined by
\beq
\Psi_0(\Bx) = \sum_{j=1}^N C_j \Big(P^{(j)}(\Bx) - \fr{\mbox{cap } F^{(j)} }{|\Bx - \BO^{(j)}|} \Big), 
\eequ{Psi0}
satisfies
$$
R(\Bx) + \Psi_0(\Bx) = 0
$$
quasi-everywhere on $\prt \GO$, which follows from \eq{as_form} and \eq{f2},  \eq{f2u}.

We set $B_\rho^{(k)} = \{ \Bx: ~|\Bx - \BO^{(k)}| < \rho \}$, and define the capacitary potential of  $F^{(k)}$ relative to $B^{(k)}_{d/4}$,
that is a unique variational solution of the Dirichlet problem
\begin{eqnarray}
\GD \tilde{P}_k (\Bx) &&= 0, ~~ \Bx \in B_{d/4}^{(k)} \setminus F^{(k)}, \l{Pt1} \\
\tilde{P}_k(\Bx) &&= 1, ~~ \Bx \in \prt ({\Bbb R}^3 \setminus F^{(k)}), \l{Pt2} \\
\tilde{P}_k (\Bx) &&= 0, ~~ |\Bx-\BO^{(k)}| = d/4. \l{Pt3}
\end{eqnarray}

%
Also, let a surface $S_d$ be a smooth perturbation of $\prt \GO$ such that
$$
S_d \subset \GO~~ \mbox{and}~~ d/4 \leq \mbox{dist} (S_d, \Bx) \leq d/2 ~~
\mbox{for all} ~~ \Bx \in \prt \GO.
$$ 
In turn, the set of all points placed between the surfaces $\prt \GO$ and $S_d$ is denoted by $\Pi_d$,  and  the function $\tilde{P}_0$ is defined as a unique variational solution of the Dirichlet problem
\begin{eqnarray}
\GD \tilde{P}_0 (\Bx) &&= 0, ~~ \Bx \in \Pi_d, \l{Pt1a} \\
\tilde{P}_0(\Bx) &&= 1, ~~ \Bx \in \prt \GO, \l{Pt2b} \\
\tilde{P}_0 (\Bx) &&= 0, ~~ \Bx \in S_d. \l{Pt3c}
\end{eqnarray}

We note that
\beq
R(\Bx) + \sum_{k=0}^N
 \tilde{P}_k (\Bx) \Psi_k(\Bx) 
 \eequ{R_Psi}
 vanishes quasi-everywhere on $\prt \GO_N$ and that the Dirichlet integral of \eq{R_Psi}  over $\GO_N$ is finite. Therefore, by harmonicity of $R$
 $$
 \int_{\GO_N} \Grad R(\Bx) \cdot \Grad \Big( R(\Bx) + \sum_{0 \leq k \leq N} \tilde{P}_k(\Bx) \Psi_k(\Bx)  \Big) d \Bx = 0.
 $$
 Hence
 $$
 \| \Grad R \|^2_{L_2(\GO_N)} \leq   \| \Grad R \|_{L_2(\GO_N)}~ \|\Grad \sum_{0 \leq k \leq N} \tilde{P}_k 
\Psi_k \|_{L_2(\GO_N)},
$$
which is equivalent to the estimate
%
 \beq
 ~~~~~~~~~~~~~~~\| ~\Grad R
  ~\|_{L_2 (\GO_N)} \leq  
 \Bigg( \sum_{k=1}^N \|~  \Grad (\tilde{P}_k \Psi_k)  ~ \|^2_{L_2 (B_{d/4}^{(k)} )} 
 + \|~  \Grad (\tilde{P}_0 \Psi_0)  ~ \|^2_{L_2 (\Pi_{d} )}  \Bigg)^{1/2}. 
 \eequ{grad_R_N_1}
 
 
 In the remaining part of the proof, we 
 obtain an upper estimate for the right-hand side in \eq{grad_R_N_1}.
 
 The inequality \eq{grad_R_N_1} and the definition of $\Psi_k$ lead to 
 $$
 \| ~\Grad R
  ~\|_{L_2 (\GO_N)}^2 \leq 2 \Big(   {\cal K}^{(1)} +    {\cal K}^{(2)} +   {\cal L}^{(1)}  +   {\cal L}^{(2)}  +   {\cal M}^{(1)}  +   {\cal M}^{(2)} + \CN   + {\cal Q} \Big),
 $$
 where
 \begin{eqnarray}
 {\cal K}^{(1)} &&= \sum_{k=1}^N \|    \Grad \Big(  \tilde{P}_k  ( v_f(\cdot) - v_f(\BO^{(k)})     ) \Big)     \|^2_{L_2 (B_{3 \Gve} ^{(k)})},
 \l{frk_K} \\
  {\cal L}^{(1)} &&= \sum_{k=1}^N \|  \sum_{1 \leq j \leq N, ~j \neq k} C_j \Grad \Big( \tilde{P}_k ( P^{(j)} (\cdot)  - 
  \fr{
  \mbox{cap} (F^{(j)}) }{|  \BO^{(k)}  - \BO^{(j)}  |}     )         \Big) \|^2_{L_2 (B_{3 \Gve} ^{(k)})},
   \l{frk_L} \\
     {\cal M}^{(1)} &&= (4 \pi)^2 \sum_{k=1}^N \Big\|  \sum_{j=1}^N C_j ~\mbox{cap} (F^{(j)})~ \Grad \Big( \tilde{P}_k ( H(\cdot, \BO^{(j)})  \nonumber \\
     &&~~~~~~- H(\BO^{(k)}, \BO^{(j)})    ) \Big)  \Big\|^2_{L_2 (B_{3 \Gve} ^{(k)})},
   \l{frk_M} \\
        {\cal N} &&= (4 \pi)^2 \sum_{k=1}^N |C_k|^2 ~\Big(\mbox{cap} (F^{(k)}) \Big)^2~ \Big(H (\BO^{(k)},  \BO^{(k)})\Big)^2  \big\| \Grad  \tilde{P}_k  \big\|^2_{L_2 (B_{d/4} ^{(k)})},  \l{frk_N} \\
        {\cal Q} && =  \|  \sum_{1 \leq j \leq N} C_j \Grad \Big( \tilde{P}_0 ( P^{(j)} (\cdot)  - 
  \fr{
  \mbox{cap} (F^{(j)}) }{|  \Bx  - \BO^{(j)}  |}     )         \Big) \|^2_{L_2 (\Pi_d)},
\l{frk_Q}
 \end{eqnarray}
 and $ {\cal K}^{(2)}, ~  {\cal L}^{(2)}, ~ \CM^{(2)}$ are defined by replacing $B_{3 \Gve}^{(k)}$ in the definitions of  $ {\cal K}^{(1)}, ~  {\cal L}^{(1)},
 ~\CM^{(1)}$ by $B_{d/4}^{(k)} \setminus B_{3 \Gve}^{(k)}$.    
 
 We start with the 
 sum $ {\cal K}^{(1)}$. Clearly,
 \begin{eqnarray}
  {\cal K}^{(1)}  &&\leq C \| \Grad v_f \|^2_{L_\infty(\Go
  )} ~\sum_{k=1}^N \int_{B_{3 \Gve}^{(k)}} \Big\{  |\Grad \tilde{P}_k(\Bx) |^2 ~|\Bx - \BO^{(k)}|^2 + \Big(  \tilde{P}_k(\Bx) \Big)^2       \Big\} d \Bx
 \nonumber \\
 && \leq C \| \Grad v_f \|^2_{L_\infty(\Go
 )} ~\sum_{k=1}^N \Gve^2 ~\mbox{cap} (F^{(k)}) 
 \l{K_sq}
 \end{eqnarray}
 and hence
 \beq
   {\cal K}^{(1)}  \leq C 
   {\Gve^3}{d^{-3}} \| \Grad v_f \|^2_{L_\infty (\Go
 )}.
 \eequ{K_sq_a}
 
 Furthermore, by Green's formula and by \eq{f1u} we deduce
 \begin{eqnarray}
   {\cal K}^{(2)}  = - \sum_{k=1}^N&& \int_{B_{d/4}^{(k)} \setminus B_{3 \Gve}^{(k)}} \tilde{P}_k(\Bx) \Big(   v_f(\Bx)  - v_f(\BO^{(k)})\Big) \Big\{   - \tilde{P}_k(\Bx) f(\Bx)  \nonumber \\
   &&+ 2 \Grad \tilde{P}_k(\Bx) \cdot 
   \Grad v_f(\Bx)   \Big\} d \Bx \nonumber \\
   - \sum_{k=1}^N  &&
 \int_{ \prt B^{(k)}_{3 \Gve}  }  \tilde{P}_k(\Bx) \Big(   v_f(\Bx)  - v_f(\BO^{(k)})\Big)  \Big\{  \tilde{P}_k(\Bx) \fr{\prt v_f}{\prt |\Bx|} (\Bx)  \nonumber \\
 && +
 (v_f(\Bx) - v_f(\BO^{(k)})) \fr{\prt \tilde{P}_k}{\prt |\Bx|} (\Bx) \Big\} d S \l{frk_K*}
 \end{eqnarray}

  By the mean value theorem for harmonic functions and the inequality $\tilde{P}_k(\Bx) \leq P^{(k)}(\Bx)$, we have
 $$
 |\Grad \tilde{P}_k(\Bx)| \leq \fr{C}{|  \Bx - \BO^{(k)}  |} ~ \max_{\By \in B} 
 {P}^{(k)}(\By),
 $$
 where $B = \{ \By: ~ |\By - \Bx| < |\Bx - \BO^{(k)}| /4    \}.
 $
 Making use of the asymptotics \eq{P} far from $\BO^{(k)}$ we deduce
%
 \beq
 |\Grad \tilde{P}_k(\Bx) | \leq C ~ \fr{\mbox{cap} (F^{(k)}) }{|\Bx - \BO^{(k)}|^2}, ~ \Bx \in B_{d/4}^{(k)} \setminus B_{3 \Gve}^{(k)}.
 \eequ{VA}
 
 
 Now we turn to the estimate of \eq{frk_K*}.  The volume integral in the right-hand side of  \eq{frk_K*} does not exceed
 \begin{eqnarray}
  C && \sum_{k=1}^N
 \| \Grad v_f \|_{L_\infty(\Go 
 )} \mbox{cap} (F^{(k)}) \int_{B_{d/4}^{(k)} \setminus B_{3 \Gve}^{(k)}} 
 \Big\{
 \fr{\mbox{cap} ( F^{(k)})}{|  \Bx - \BO^{(k)}|} |f(\Bx)| \nonumber \\
 && + \fr{\mbox{cap} ( F^{(k)})}{|\Bx - \BO^{(k)}|^2} \| \Grad v_f  \|_{L_\infty(\Go 
 ) } \Big\} d \Bx \nonumber \\
 && \leq  C \Gve d^{-3}  \| \Grad v_f \|_{L_\infty(\Go 
 )} \Big\{ \Gve d^2 \|f\|_{L_\infty(\GO_N) }+ \Gve d \| \Grad v_f \|_{L_\infty(\Go 
 ) }       \Big\} \nonumber \\
 && \leq C  {\Gve^2}{d^{-2}}  \| f \|^2 _{L_\infty(
 \GO_N
 )}. \l{Grad_1}
 \end{eqnarray}
 
 By  $0 \leq \tilde{P}_k(\Bx) \leq 1$ and \eq{VA}, the surface integral in  \eq{frk_K*} is dominated by
 \beq
 C \Gve^3 d^{-3} \| \Grad v_f \|_{L_\infty (\Go 
 )}^2 \leq C \Gve^3 d^{-3} \|f\|^2_{L_\infty(\GO_N)}.
 \eequ{Grad_2}
 Combining \eq{Grad_1} and \eq{Grad_2} we arrive at the estimate
 \beq
 {\cal K}^{(2)} \leq C \Gve^2 d^{-2} \|f \|^2_{L_\infty(\GO_N)}.
 \eequ{frk_K_grad}
 
 Let us estimate ${\cal L}^{(1)}$ (see \eq{frk_L}). Obviously,
 $$
 {\cal L}^{(1)} \leq \sum_{k=1}^N \Bigg(  \sum_{1 \leq j \leq N, ~ j \neq k} |C_j | ~    
 \Big\| \Grad  \Big(   \tilde{P}_k  (  P^{(j)}(\cdot) - 
 \fr{
 \mbox{cap} (F^{(j)})}{|  \BO^{(k)} - \BO^{(j)} |}
  )  \Big)     \Big\|_{L_2(B_{3 \Gve}^{(k)})}      \Bigg)^2.
  $$
  Furthermore, when $j \neq k$ we have
 \beqa
 && 
  \Big\| \Grad  \Big(   \tilde{P}_k  (  P^{(j)}(\cdot) - \fr{
  \mbox{cap} (F^{(j)})}{|  \BO^{(k)} - \BO^{(j)} |}  
   )  \Big)     \Big\|_{L_2(B_{3 \Gve}^{(k)})} 
 \nonumber \\
  && \leq  \| (  \Grad \tilde{P}_k  )   \Big(  P^{(j)}(\cdot )  - \fr{
  \mbox{cap} (F^{(j)}) }{|  \BO^{(k)} - \BO^{(j)} |}   
   \Big)   \|_{L_2(B_{3 \Gve}^{(k)})}
 + \|   \tilde{P}_k   \Grad  P^{(j)}  \|_{L_2(B_{3 \Gve}^{(k)})},
 \nonumber
 \end{eqnarray}
 which does not exceed
 $$
 C \Bigg\{     \fr{\Gve ~\mbox{cap} (F^{(j)})  (\mbox{cap} (F^{(k)}))^{1/2}    }{|  \BO^{(k)} - \BO^{(j)}   |^2}    +
 \fr{ \Gve^{3/2} \mbox{cap} (F^{(j)})}{|  \BO^{(k)} - \BO^{(j)}  |^2}   \Bigg\}.
 $$
 Hence, using Lemma 1 we deduce
 \beqa
 && {\cal L}^{(1)} \leq  C~\sum_{k=1}^N \Bigg(   \Gve^2 \sum_{1 \leq j \leq N, ~ j \neq k} |C_j| \fr{(  \mbox{cap} (F^{(j)})   )^{1/2}}{ |   \BO^{(k)} - \BO^{(j)}  |^2}    \Bigg)^2
 \nonumber \\
 && \leq C ~\Gve^4 \sum_{k=1}^N
 \Big\{    \sum_{1 \leq j \leq N, ~ j \neq k} C_j^2 ~ \mbox{cap} ( F^{(j)})   \sum_{1 \leq j \leq N, ~ j \neq k} \fr{1}{|   \BO^{(k)} - \BO^{(j)}  |^4}      \Big\}, 
 \nonumber
 \end{eqnarray}
 and therefore
 \beq
  {\cal L}^{(1)} \leq C \Gve^{5} d^{-10} \| v_f \|_{L_\infty(\Go 
 )}^2.
 \eequ{star_sq}
 
 Similar steps can be followed to estimate $\CQ$ in \eq{frk_Q}.  We have
\beqa
&&  \Big\| \Grad  \Big(   \tilde{P}_0  (  P^{(j)}(\cdot) - \fr{
  \mbox{cap} (F^{(j)})}{|  \Bx - \BO^{(j)} |}  
   )  \Big)     \Big\|_{L_2(\Pi_d)} 
 \nonumber \\
 &&  \leq  \Big\| (  \Grad \tilde{P}_0  )   \Big(  P^{(j)}(\cdot )  - \fr{
  \mbox{cap} (F^{(j)}) }{|  \Bx - \BO^{(j)} |}   
   \Big)   \Big\|_{L_2(\Pi_d)}
 + \Big\|   \tilde{P}_0   \Big(\Grad  P^{(j)} +   \mbox{cap}(F^{(j)})  \fr{\Bx-\BO^{(j)}}{|\Bx-\BO^{(j)}|^3}\Big)  \Big\|_{L_2(\Pi_d)},
 \nonumber
 \end{eqnarray}
 which does not exceed
$$
 C  {\Gve ~\mbox{cap} (F^{(j)})    }{|\Bx - \BO^{(j)}|^{-2} }  , ~~ \Bx \in \Pi_d.
 $$
 As above, we use Lemma 1 to deduce 
 \beqa
&& {\cal Q} \leq  C~\Gve^3 \Bigg(   \sum_{1 \leq j \leq N} |C_j| \fr{(  \mbox{cap} (F^{(j)})   )^{1/2}}{ |   \Bx - \BO^{(j)}  |^2}    \Bigg)^2
 \nonumber \\
 && \leq C ~\Gve^3
  \sum_{1 \leq j \leq N} C_j^2 ~ \mbox{cap} ( F^{(j)})   \sum_{1 \leq j \leq N} \fr{1}{|   \Bx - \BO^{(j)}  |^4}  \leq   C \Gve^{4} d^{-7} \| v_f \|_{L_\infty(\Go 
 )}^2 ~~ \mbox{for}~~ \Bx \in \Pi_d.
\l{Q_est}
\end{eqnarray}

Next, we estimate ${\cal L}^{(2)}$. 
Integration by parts gives 
\beqa
&& \hspace{-.5cm} \int_{  B_{d/4}^{(k)} \setminus B_{3 \Gve}^{(k)}    } \Grad \Big(   \tilde{P}_k(\Bx)
(  P^{(j)}(\Bx)  -       \fr{
\mbox{cap} (F^{(j)}) }{|  \BO^{(k)} - \BO^{(j)} |}        
)
 \Big) 
\cdot \Grad    \Big(   \tilde{P}_k(\Bx)
(  P^{(m)}(\Bx)  -  \fr{
\mbox{cap} (F^{(m)})  }{|  \BO^{(k)} - \BO^{(m)} |}              
  ) \Big)d \Bx
\nonumber \\
&& =- 2 \int_{   B_{d/4}^{(k)} \setminus B_{3 \Gve}^{(k)}   }
\tilde{P}_k(\Bx) \Big(  P^{(j)} (\Bx) -      \fr{
\mbox{cap} (F^{(j)} ) }{|  \BO^{(k)} - \BO^{(j)} |}        
 \Big) \Big(   \Grad \tilde{P}_k(\Bx ) \cdot \Grad P^{(m)}(\Bx)  \Big) d \Bx 
\nonumber \\
&& -\int_{\prt B_{3 \Gve}^{(k)}}
 \tilde{P}_k(\Bx)~\Big(  P^{(j)}(\Bx)  -  \fr{
 \mbox{cap} (F^{(j)}) }{|  \BO^{(k)} - \BO^{(j)} |}  
 \Big)
\Big\{  \tilde{P}_k(\Bx)~ \fr{\prt}{\prt   |\Bx|} 
P^{(m)}(\Bx) 
\nonumber \\
&& + ( P^{(m)}(\Bx)  - 
\fr{
\mbox{cap} (F^{(m)}) }{|  \BO^{(k)} - \BO^{(m)} |}   ) \fr{\prt \tilde{P}_k}{\prt  |\Bx|} 
(\Bx)       \Big\} dS,
\l{ring_00}
\end{eqnarray}

When $j \neq k$ and $m \neq k$, the volume integral in the right-hand side of \eq{ring_00} is estimated as follows
\beqa
&& \Big| \int_{   B_{d/4}^{(k) } \setminus   B_{3 \Gve}^{(k) }    }
\tilde{P}_k(\Bx) \Big(  P^{(j)}(\Bx) - 
\fr{
\mbox{cap} (F^{(j)}) }{|  \BO^{(k)} - \BO^{(j)} |}  \Big) \Big(   \Grad \tilde{P}_k(\Bx ) \cdot \Grad P^{(m)}(\Bx)  \Big) d \Bx  \Big|
\nonumber \\
&& \leq C 
~ \fr{\mbox{cap} (F^{(k)}) ~ \mbox{cap} (F^{(j)})  }{|  \BO^{(k)} - \BO^{(j)}|^2}  ~ \int_{    B_{d/4}^{(k) } \setminus   B_{3 \Gve}^{(k) }       } \fr{\mbox{cap} (F^{(k)}) }{|\Bx- \BO^{(k)}|^2}  
~\fr{\mbox{cap} ( F^{(m)}) }{|\Bx- \BO^{(m)}|^2} d\Bx 
\l{volume_ring} \\
&& \leq C
\fr{(\mbox{cap} (F^{(k)}) )^2 \mbox{cap} (F^{(j)}) \mbox{cap} (F^{(m)})  d}{|  \BO^{(k)} - \BO^{(j)}|^2 ~|  \BO^{(k)} - \BO^{(m)}|^2} 
\leq C 
\fr{\Gve^2 d ~ \mbox{cap} ( F^{(j)})~ \mbox{cap} (F^{(m)})  }{|  \BO^{(k)} - \BO^{(j)}|^2 ~|  \BO^{(k)} - \BO^{(m)}|^2} 
\nonumber
\end{eqnarray}

In turn, when $j \neq k$ and $m \neq k$ the modulus of the surface integral in the right-hand side of \eq{ring_00} does not exceed
\beqa
&& C
\fr{\Gve ~\mbox{cap} (F^{(j)}) }{|\BO^{(k)} - \BO^{(j)}|^2} \int_{ \prt B_{3 \Gve}^{(k)}
 }
\Bigg\{    
\fr{\mbox{cap} (F^{(m)}) }{| \Bx - \BO^{(m)}   |^2} + \fr{\Gve ~\mbox{cap} (F^{(m)})  }{|  \BO^{(m)}   - \BO^{(k)}    |^2} \Big|  \fr{\prt \tilde{P}_k}{\prt |\Bx| } 
(\Bx)  \Big|
\Bigg\} dS
\nonumber \\
&& \leq C 
\fr{\Gve ~\mbox{cap} ( F^{(j)}) ~\mbox{cap} (F^{(m)})  }{|\BO^{(k)} - \BO^{(j)}|^2 |\BO^{(k)} - \BO^{(m)}|^2}  
 \Bigg\{     \Gve^2 + \Gve \int_{  \prt B_{3 \Gve}^{(k) } } 
 \fr{\mbox{cap} (F^{(k)}) }{|\Bx - \BO^{(k)}|^2}   dS       \Bigg\}
\nonumber \\
&& \leq C 
\fr{\Gve^3 ~\mbox{cap} (F^{(j)}) ~\mbox{cap} (F^{(m)}) }{|\BO^{(k)} - \BO^{(j)}|^2 |\BO^{(k)} - \BO^{(m)}|^2}  .
\l{surface_ring}
\end{eqnarray}

We have
\beqa
 {\cal L}^{(2)} = &&\sum_{1 \leq m,j \leq N} C_m C_j
\sum_{1 \leq k \leq N, ~ k \neq m, k \neq j} 
\int_{B_{d/4}^{(k)}\setminus B_{3 \Gve}^{(k)}} 
\Grad \Big(   \tilde{P}_k(\Bx)  (  P^{(j)}(\Bx) -      \fr{
\mbox{cap} (F^{(j)}) }{|  \BO^{(k)} - \BO^{(j)} |}           
 )  \Big) 
\nonumber \\
&&  \cdot \Grad \Big(  \tilde{P}_k(\Bx)  (  P^{(m)}(\Bx) -         \fr{
 \mbox{cap} (F^{(m)}) }{|  \BO^{(k)} - \BO^{(m)} |}           
  )     \Big)   d \Bx,
\nonumber
\end{eqnarray}
and 
by  \eq{ring_00} ,  \eq{volume_ring} and \eq{surface_ring}
\beqa
&& \hspace{-.5cm}{\cal L}^{(2)} \leq C \Gve^2 d \sum_{1 \leq m,j \leq N} |C_m |   |C_j|  \sum_{1 \leq k \leq N, ~ k \neq m, k \neq j} \fr{\mbox{cap} (F^{(j)})  \mbox{cap} ( F^{(m)}) }{|  \BO^{(k)} - \BO^{(j)}  |^2 ~ |  \BO^{(k)} - \BO^{(m)}    |^2}
\nonumber \\
&&  \hspace{-.5cm} = C ~ \fr{\Gve^2}{d^2}  {\small \sum_{1 \leq m,j \leq N} |C_m |   |C_j| }
{\small \mbox{cap} (F^{(j)})  ~ \mbox{cap} (F^{(m)})
\sum_{\scriptsize 
1 \leq k \leq N, ~
 k \neq m, k \neq j 
}  \fr{d^3}{ \small |  \BO^{(k)} - \BO^{(j)}  |^2 ~ |  \BO^{(k)} - \BO^{(m)}    |^2}   },
\nonumber
\end{eqnarray}
and therefore
\beq
\CL^{(2)} \leq C ~ \fr{\Gve^2}{d^2}  ~ \sum_{1 \leq m,j \leq N} \fr{ |C_m |   |C_j|  \mbox{cap} ( F^{(j)}) ~  \mbox{cap} ( F^{(m)}) }{d + |  \BO^{(j)} - \BO^{(m)}  |}.
\eequ{square_star}
 Let us introduce a piece-wise constant function
$$
\xi(\Bx) = \left\{    
\begin{array}{cc}
|C_m| (\mbox{cap} ( F^{(m)} ))^{1/2} , & \mbox{when }~ \Bx \in B^{(m)}_{d/4}, \\
0, & \mbox{otherwise}.
\end{array}
\right.
$$
Then the inequality \eq{square_star} 
leads to 
\beqa
{\cal L}^{(2)} \leq && 
C ~ \fr{\Gve^3}{d^8} ~ \sum_{1 \leq m,j \leq N} \fr{  \Big(  |C_m| (\mbox{cap} (F^{(m)}) )^{1/2} \Big) ~\Big(  |C_j| (\mbox{cap} (F^{(j)}) )^{1/2} \Big)  d^6  }{   d + |\BO^{(j)}- \BO^{(m)}|     }
\nonumber \\
\leq && C   
~  \fr{\Gve^3}{d^8} ~ \int_\Go \int_\Go \fr{ \xi(\BX) \xi(\BY) }{d + |  \BX - \BY   |  } d\BX d \BY \leq 
C~ \fr{\Gve^3}{d^8} ~\|\xi\|^2_{L_2(\Go)},
\nonumber 
\end{eqnarray}
where the constant $C$ depends on $\Go$, and using Lemma 1 we deduce
\beq
{\cal L}^{(2)} \leq 
C~ \fr{\Gve^3}{d^8} \sum_{1 \leq j \leq N} C_j^2 ~\mbox{cap} (F^{(j)})  d^3 \leq 
C~\fr{\Gve^4}{d^8}   \| v_f\|^2_{L_2(\Go)}. 
\eequ{S1_result}

To evaluate $\CM^{(1)} + \CM^{(2)}$ we apply the result of Lemma 2 and use the same algorithm as for  $\CK^{(1)}$  and $\CK^{(2)}$ to deduce
\beq
\CM^{(1)} + \CM^{(2)} \leq C \|v_f\|^2_{L_\infty(\Go)} 
{\Gve}{d^{-3}} \Big(  \Gve^3 d^{-3} + \Gve^2 d^{-2}    \Big) \leq C \Gve^3 d^{-5} \|v_f\|^2_{L_\infty(\Go)}.
\eequ{CM_est}

Similarly, applying Lemma 2, we derive the estimate for the term $\CN$ 
\beq
\CN \leq C \Gve^3 d^{-3} \|v_f \|^2_{L_\infty (\Go)}.
\eequ{CN_est}

The proof is completed by the reference to   \eq{K_sq},  \eq{frk_K_grad},  \eq{star_sq},  \eq{Q_est}, \eq{S1_result} , \eq{CN_est}. 
$\Box$

 \section{Meso scale approximation of Green's function in $\GO_N$}
 \l{meso_green}
 
 Let $G_N(\Bx, \By)$ be Green's function of the Dirichlet problem for the operator $-\GD$ in $\GO_N$. In this section, we derive the 
 asymptotic approximation of $G_N(\Bx, \By)$ 
 and estimate the remainder term. In the asymptotic algorithm, we  will refer to the 
 algebraic system 
 similar to that 
 of Section \ref{aux_syst}. 
 We need here Green's functions $g^{(j)}(\Bx, \By)$ of the Dirichlet problem for the operator   $-\GD$ in ${\Bbb R}^3 \setminus F^{(j)}, ~ j=1,\ldots,N.$ The notation $h^{(j)}$ will be used for the regular part of $g^{(j)}$, that is
 \beq
 h^{(j)} (\Bx, \By) = 
 ({4 \pi |  \Bx - \By | })^{-1}    - g^{(j)}( \Bx,  \By  ), ~~ \Bx, \By \in {\Bbb R}^3 \setminus F^{(j)}.
 \eequ{h_reg_part}
According to Lemma 2 of \cite{JCAM}, 
the functions $h^{(j)}$ allow for the following estimate:
\beq
\Big| h^{(j)}(\Bx, \By) - \fr{P^{(j)}(\By)}{4 \pi  |\Bx - \BO^{(j)}|}  \Big| \leq ~\mbox{const}~ \fr{\Gve P^{(j)}(\By)}{|\Bx - \BO^{(j)}|^2},
\eequ{h_asymp}
for all $\By \in {\Bbb R}^3 \setminus F^{(j)}$ and $|\Bx - \BO^{(j)}| > 2 \Gve.$
 
 The principal result of this section is
 
 {\bf Theorem 4.}  {\it Let the small parameters $\Gve $ and $d$, introduced in {\rm \eq{eps_d}}, satisfy the inequality $\Gve < c ~ d^2,$
where $c$ is a sufficiently small absolute constant.
Then  
 \beqa
&&  G_N(\Bx, \By) = G(\Bx, \By) - \sum_{j=1}^N \Bigg\{ h^{(j)} (\Bx, \By) 
  \l{G_N_as}  \\
 && - P^{(j)}(\By) H(\Bx, \BO^{(j)}) 
  - P^{(j)}(\Bx) H(\BO^{(j)}, \By)      
  +4 \pi~\mbox{\rm cap}(F^{(j)}) H(\Bx, \BO^{(j)}) H(\BO^{(j)}, \By)    
\nonumber \\
&& +H(\BO^{(j)}, \BO^{(j)}) ~ T^{(j)} (\Bx )  T^{(j)} (\By)   
-\sum_{i=1 }
^N 
 {\cal C}_{ij}  T^{(i)} (\Bx )  T^{(j)} (\By)   \Bigg\} + \CR(\Bx, \By),
\nonumber
\end{eqnarray}
where
\beq
T^{(j)}(\By) =  P^{(j)}(\By) -    4 \pi ~\mbox{\rm cap}(F^{(j)})  H(\BO^{(j)}, \By),
\eequ{T_def}
with the capacitary potentials $P^{(j)}$ and the regular part $H$ of Green's function $G$  of $\GO$ being  the same as in Section {\rm \ref{aux_pr}}. The matrix  $ \BCC = ({\cal C}_{ij})_{i,j=1}^N$ 
is defined by
\beq
\BCC = (\BI + \BS \BD)^{-1} \BS,
\eequ{CC_def}
where 
$\BS$ and $\BD$ are the same as in {\rm \eq{S_matrix1}, \eq{D_matrix1}}. The remainder $\CR(\Bx, \By)$ is a harmonic function, both in $\Bx$ and $\By$, 
and satisfies the estimate
\beq
|\CR(\Bx, \By)| \leq 
{\mbox{\rm const}~ \Gve d^{-2}} 
\eequ{CR_est}
uniformly with respect to $\Bx$ and $\By$ in $\GO_N$. 
 }
 
 \vspace{.2in}
 
 Prior to the proof of the theorem, we formulate an auxiliary result.
 
  \vspace{.2in}
  
 {\bf Lemma 3.} {\it Let the small parameters $\Gve$ and $d$, defined in   {\rm \eq{eps_d}},  obey the inequality {\rm \eq{d2_constr}}. Then the matrix
$\BCC$ in 
{\rm \eq{CC_def}} satisfies the estimate
 \beq
 \| \BCC \|_{{\Bbb R}^N \to {\Bbb R}^N} \leq c  d^{-3}, 
 \eequ{l3_eq1}
 where $c$ is an absolute constant.}
 
 {\it Proof.} 
 First, we note that
 \beq
  \| \BCC \|_{{\Bbb R}^N  \to {\Bbb R}^N} \leq \mbox{const }  \| \BS \|_{{\Bbb R}^N \to {\Bbb R}^N} ,
 \eequ{l3_eq2}
 which follows from Lemma 2, where $\BV_f$ should be replaced by the columns of the matrix $\BS$. 
 
 Additionally, 
 \beq
 \|\BS \|_{{\Bbb R}^N  \to {\Bbb R}^N} \leq \mbox{const }  d^{-3}.
 \eequ{T4_eq_a}
 To  verify this estimate we introduce a vector $\BGx = (\xi_j)_{j=1}^N, ~ \|\BGx\| = 1,$ and a function $\xi(\Bx)$ defined in $\GO$  by
\beq
\xi(\Bx) =    \left\{ \begin{array}{cc}
\xi_j & ~\mbox{in} 
~ B^{(j)} = \{\Bx: |\Bx - \BO^{(j)}| < d \} ,  ~ j = 1, \ldots, N,\\
0 & ~ \mbox{otherwise.}
\end{array} \right.
\eequ{T4_eq_b}
Then
\beqa
&& \langle \BS \BGx, \BGx \rangle \leq \mbox{const} ~ d^{-6}  \int_\Go \int_\Go \fr{\xi (\BX)  \xi(\BY)  d \BX d \BY}{4 \pi |\BX - \BY|}
\leq \mbox{const} ~ d^{-6} \int_\Go |\xi (\BX)|^2 d X \leq \mbox{const}~ d^{-3},
\nonumber
\end{eqnarray}
which yields   \eq{T4_eq_a}. Then  \eq{l3_eq2}  together with  \eq{T4_eq_a} lead to  \eq{l3_eq1}. $\Box$
  \vspace{.2in}
 
 {\it Proof of Theorem 4.}   The harmonicity of $\CR$ follows directly from \eq{G_N_as}.


 Let us 
 estimate the boundary values of $\CR$ on $\prt \GO_N$.
 
 If $\Bx \in \prt \GO$ and $\By \in \GO_N$, then according to  the definitions of Section \ref{Green_unpert} for Green's function of $\GO$ and its regular part
 the remainder term $\CR$ in \eq{G_N_as} takes the form
 \beqa
  \CR(\Bx, \By) = && \sum_{j=1}^N \Bigg\{ h^{(j)} (\Bx, \By) 
 - 
 \fr{ P^{(j)}(\By)}{4 \pi |\Bx - \BO^{(j)}|}   
  - H(\BO^{(j)}, \By) \Big( P^{(j)}(\Bx)       
  -
  \fr{\mbox{\rm cap}(F^{(j)})}{|\Bx - \BO^{(j)}|} 
   \Big)
\nonumber \\
&& +H(\BO^{(j)}, \BO^{(j)}) ~ 
T^{(j)} (\By)  \Big(   P^{(j)}(\Bx) -    
 ~\fr{\mbox{\rm cap}(F^{(j)})}{|  \Bx - \BO^{(j)}    |}  
 \Big)  
\nonumber \\
&& -\sum_{i=1} 
^N  {\cal C}_{ij}  
T^{(j)} (\By)     \Big(   P^{(i)}(\Bx) -    
 ~\fr{\mbox{\rm cap}(F^{(i)})}{|  \Bx - \BO^{(i)}    |}  
 \Big)      \Bigg\} .
 \nonumber
 \end{eqnarray}
 Taking into account the estimate \eq{h_asymp} 
 for $h^{(j)}$ together with the asymptotic representation \eq{P} of $P^{(j)}$ 
 we obtain
 %
 \beq
   \CR(\Bx, \By) = \sum_{j=1}^N O\Big(\fr{\Gve P^{(j)}(\By)}{  | \Bx - \BO^{(j)}  |^2 }\Big)  + \sum_{j=1}^N  \sum_{i=1}^N  {\cal C}_{ij} T^{(j)}(\By)
   O \Big( 
   \fr{\Gve^2 
   }{ 
   |  \Bx - \BO^{(i)}  |^2}    \Big) 
 \eequ{T3_eq1}
 for all $\Bx \in \prt \GO$. 
 
 Here  $\CR(\Bx, \By)$ is harmonic as a function of $\By$. Next, we estimate \eq{T3_eq1} for $\By \in \prt \GO_N$.
 
 If $\By, \Bx \in \prt \GO$ then  \eq{T3_eq1}, \eq{T_def} and \eq{P} lead to
  \beq
   \CR(\Bx, \By) = \sum_{j=1}^N O\Big(\fr{\Gve^2}{ | \By - \BO^{(j)}  |  | \Bx - \BO^{(j)}  |^2 }\Big) + \sum_{j=1}^N  \sum_{i=1}^N  {\cal C}_{ij} 
   O \Big( 
   \fr{\Gve^4 
   }{ 
   |  \Bx - \BO^{(i)}  |^2  |  \By - \BO^{(j)}  |^2}    \Big) .
    \eequ{T3_eq1_a}
    
  Using  \eq{l3_eq1} we can estimate the double sum from  \eq{T3_eq1_a}. For a fixed $\Bx \in \prt \GO$, let us introduce a vector
  $$
  \BV = \Big(   \fr{\Gve^2}{|\Bx - \BO^{(i)}|^2 }   \Big)_{i=1}^N,
  $$
  and a function $V(\BX)$ defined in $\GO$ by
  $$
  V(\BX) = \left\{ \begin{array}{cc}
  V_j & \mbox{when} ~ \BX \in  B^{(j)}, \\
  0 & \mbox{otherwise},
  \end{array}
  \right.
  $$
  where the balls $B^{(j)}$ are the same as in \eq{T4_eq_b}. It follows from Lemma 3 that the double sum in \eq{T3_eq1_a}  does not exceed
  $$
  c ~d^{-3} \|\BV \|^2 \leq \fr{\mbox{const }}{d^6} \int_\Go (V(\BX))^2 d \BX \leq \fr{\mbox{const } \Gve^4}{d^7}.
  $$
  The above estimate together with \eq{T3_eq1_a} imply
  \beq
  \CR(\Bx, \By) = O(\Gve^2 d^{-3} |\log d| + \Gve^4 d^{-7}) ~\mbox{when}~ \Bx, \By \in \prt \GO.
  \eequ{T4_eq_d}
  
  Now, we estimate \eq{T3_eq1} for $\By \in \prt ( {\Bbb R}^N \setminus F^{(m)})$. In this case we have
  \beqa
   \CR(\Bx, \By) =&& O\Big(\fr{\Gve}{  | \Bx - \BO^{(j)}  |^2 }\Big)   + \sum_{1 \leq j \leq N, j \neq m}  O\Big(\fr{\Gve^2}{  | \By - \BO^{(j)}  | | \Bx - \BO^{(j)}  |^2 }\Big)  
   \nonumber \\
&&   +  \sum_{i=1}^N    O \Big( 
   \fr{\Gve^2 
   }{ 
   |  \Bx - \BO^{(i)}  |^2}    \Big)  \Bigg\{ {\cal C}_{im} \Big(1 - 4 \pi \mbox{ cap} (F^{(m)})H(\BO^{(m)}, \By)\Big) 
  \nonumber \\
  &&  +   4 \pi  \sum_{1 \leq j \leq N, j \neq m}  {\cal C}_{ij}  \Big( \mbox{ cap} (F^{(j)}) G(\BO^{(j)}, \By) + O(\fr{\Gve^2}{|\By - \BO^{(j)}|^2}) \Big) \Bigg\}.
  \l{T4_eq_e}
  \end{eqnarray}
  We also note that 
   according to  \eq{CC_def}  the coefficients ${\cal C}_{ij}$ satisfy the system of algebraic equations
 \beq
 \hspace{-1cm} (1-\Gd_{im})G(\BO^{(m)}, \BO^{(i)}) - {\cal C}_{im}
  - 4 \pi \sum_{1\leq j \leq N, ~ j \neq m}   {\cal C}_{ij} 
  ~\mbox{cap} (F^{(i)}) ~ G(\BO^{(m)}, \BO^{(i)}) =0, ~~ m,i = 1,\ldots, N.
 \eequ{T3_sys}
  Hence, in the above formula \eq{T4_eq_e}  the expression in curly brackets can be written as
  \beqa
   && \hspace{-1cm}{\cal C}_{im} + O(|{\cal C}_{im}| \Gve)
   +   4 \pi  \sum_{1 \leq j \leq N, j \neq m}  {\cal C}_{ij}  \Big( \mbox{ cap} (F^{(j)}) G(\BO^{(j)}, \By) + O(\fr{\Gve^2}{|\By - \BO^{(j)}|^2}) \Big) 
   \l{T4_eq_f} \\
  &&\hspace{-1cm}  = (1-\Gd_{im}) G(\BO^{(m)}, \BO^{(i)}) + O (\Gve d^{-1}) + \sum_{1 \leq j \leq N, j \neq m}  {\cal C}_{ij} O(\fr{\Gve^2}{|\By - \BO^{(j)}|^2}), ~~ \By \in \prt ( {\Bbb R}^N \setminus F^{(m)}),
\nonumber 
  \end{eqnarray}
  and then formulae  \eq{T4_eq_e} and  \eq{T4_eq_f} imply
   $$
   \CR(\Bx, \By) = O\Big(   \Gve d^{-2}  + \Gve^2 |\log d| d^{-3} + \Gve^3 d^{-4}   \Big)   + \sum_{1 \leq i \leq N}  \sum_{1 \leq j \leq N, j \neq m} {\cal C}_{ij}  O\Big(\fr{\Gve^4}{  | \By - \BO^{(j)}  |^2 | \Bx - \BO^{(i)}  |^2 }\Big) ,
   $$
   where the estimate of the double sum is similar to  \eq{T3_eq1_a}. Thus, we obtain
   \beq
    \CR(\Bx, \By) = O\Big(   \Gve d^{-2}  + \Gve^2 |\log d| d^{-3} + \Gve^3 d^{-4}   + \Gve^4 d^{-7}\Big) = O(\Gve d^{-2}),
    \eequ{T4_eq_g} 
for all $\Bx \in \prt \GO$ and $\By \in \prt ( {\Bbb R}^3 \setminus F^{(m)} ), ~ m= 1,\ldots, N.$
 
 Using the estimates    \eq{T4_eq_d}  and    \eq{T4_eq_g} and applying the maximum principle for harmonic functions we deduce that
 \beq
   \CR(\Bx, \By) = O(   \Gve d^{-2} ),
      \eequ{T4_eq_h} 
   for all $\Bx \in \prt \GO$ and $\By \in \GO_N$.

 In turn, when $\Bx \in \prt ( {\Bbb R}^3 \setminus F^{(k)}  )$,  the 
formula \eq{G_N_as} and the definition \eq{h_reg_part} of $h^{(j)}$   
 lead to the expression for the remainder term on the boundary of the inclusion
 \beqa
  \CR(\Bx, \By) = && H(\Bx, \By) - H(\BO^{(k)}, \By) + \sum_{1 \leq j \leq N, ~ j \neq k}
  \Big( h^{(j)}(\Bx, \By) - P^{(j)}(\Bx) H(\BO^{(j)}, \By) \Big)
 \nonumber \\
 && + \sum_{j=1}^N T^{(j)} (\By) \Bigg(     H(\BO^{(j)}, \BO^{(j)}) ~ T^{(j)} (\Bx )  - H(\Bx, \BO^{(j)})   
 -\sum_{j=1}
 ^N  {\cal C}_{ij}  T^{(i)} (\Bx )   \Bigg).
 \l{T3_add1}
 \end{eqnarray}
 Using  the formulae   \eq{P}   and   \eq{h_asymp}  
 for $P^{(j)}$   and $h^{(j)}$  together with the definition \eq{T_def} 
 of
 $T^{(j)}$ and the definition of Section \ref{Green_unpert} of the regular part of Green's function of $\GO$ we deduce that 
 \beq
 h^{(j)}(\Bx, \By) - P^{(j)}(\Bx) H(\BO^{(j)}, \By) = \fr{T^{(j)}(\By)}{4 \pi |  \Bx - \BO^{(j)}   |} + O\Big(\fr{\Gve^2 + \Gve P^{(j)}(\By)}{|\BO^{(k)} - \BO^{(j)}|^2}\Big),~~ j \neq k,
 \eequ{T3_add2}
 and
 \beq
 H(\Bx, \By) = H(\BO^{(k)}, \By) + O(\Gve),
 \eequ{T3_add3}
 for $\Bx  \in \prt ( {\Bbb R}^3 \setminus F^{(k)}  )$ and $\By \in \GO_N$.
   The representations \eq{T_def} together with \eq{T3_add1}--\eq{T3_add3} 
 imply 
  \beqa
  \CR(\Bx, \By) = && 
 \sum_{1\leq j \leq N, ~ j \neq k} \Bigg\{ \fr{T^{(j)}(\By)}{ 4 \pi | \Bx - \BO^{(j)}|} + O \Big( \fr{ \Gve^2+\Gve P^{(j)}(\By)}{|  \BO^{(k)} - \BO^{(j)}|^2}  \Big) \Bigg\}
 \nonumber \\
&&  - \sum_{j=1}^N 
  T^{(j)}(\By) \Bigg(  H(\BO^{(k)},  \BO^{(j)})    -  H(\BO^{(j)}, \BO^{(j)}) T^{(j)}(\Bx) 
  \nonumber \\
 && + \sum_{i=1}^N 
  {\cal C}_{ij}  T^{(i)}(\Bx)    \Bigg)  +\sum_{j=1}^N O(\Gve |T^{(j)}(\By)|) . 
 \l{T3_add4}
 \end{eqnarray}
 Bearing in mind the asymptotic formula \eq{P} for the capacitary potentials and the definition \eq{T_def} 
 we deduce that for $\Bx \in \prt ( {\Bbb R}^3 \setminus F^{(k)}  )$
 \begin{eqnarray}
 T^{(j)} (\Bx ) &&=  
 \fr{\mbox{cap} (F^{(j)}) }{|\Bx - \BO^{(j)}|}  
 -    4 \pi ~\mbox{\rm cap}(F^{(j)})  H(\Bx, \BO^{(j)}) + O\Big(\fr{\Gve ~\mbox{cap} (F^{(j)}) }{|\Bx - \BO^{(j)}|^2}   \Big) 
 \l{T_as} \\
 && =    4 \pi ~\mbox{\rm cap}(F^{(j)})  G(\Bx, \BO^{(j)}) + O\Big(\fr{\Gve^2}{|\Bx - \BO^{(j)}|^2}   \Big)  ,    ~ j \neq k.  \nonumber
  \end{eqnarray}
 Thus, \eq{T3_add4} can be rearranged in the form
 \beqa
 \CR(\Bx, \By) =&& \hspace{-.2in} \sum_{1\leq j \leq N, ~ j \neq k} \hspace{-.2in} T^{(j)}(\By) \Bigg\{ G(\BO^{(k)}, \BO^{(j)}) - {\cal C}_{kj}
 \nonumber \\
  &&- 4 \pi \sum_{1\leq i \leq N, ~ i \neq k} 
  {\cal C}_{ij} 
  ~\mbox{cap} (F^{(i)}) ~ G(\BO^{(k)}, \BO^{(i)}) \Bigg\} + \CR^{(1)}(\Bx, \By),
 \l{T3_eq3}
 \end{eqnarray}
 where
 \beqa
 \CR^{(1)}(\Bx, \By) =&& O(\Gve) +   \sum_{1 \leq j \leq N, ~ j \neq k}  
 \Bigg\{   O\Big(   \fr{\Gve ~ |T^{(j)}(\By)|}{ |\BO^{(k)} - \BO^{(j)}|  } \Big)
 + 
 O(\Gve d^{-1} |T^{(j)}(\By)|)
 \Bigg\}
 \nonumber \\
&& + \sum_{1 \leq j \leq N, ~ j \neq k}   O\Big( \Gve |T^{(j)}(\By)| +  \fr{ \Gve^2+\Gve P^{(j)}(\By)}{|  \BO^{(k)} - \BO^{(j)}|^2}   \Big)
 \nonumber \\
 && +\sum_{1\leq i \leq N, ~ i \neq k}   
 O \Big(  \fr{\Gve ~\mbox{cap} (F^{(i)}) }{d~ |  \BO^{(k)} - \BO^{(i)}|^2} \Big)
 \l{T3_eq4} \\  
 && + \sum_{j=1}^N 
  \sum_{1\leq i \leq N, ~ i \neq k}   
   {\cal C}_{ij} T^{(j)}(\By)  
 O \Big(  \fr{\Gve ~\mbox{cap} (F^{(i)})}{|  \BO^{(k)} - \BO^{(i)}|^2} \Big)   
 \nonumber
 \end{eqnarray}
 It follows from  \eq{CC_def} that the coefficients ${\cal C}_{ij}$ satisfy the system of algebraic equations
 \beq
 (1-\Gd_{kj}) G(\BO^{(k)}, \BO^{(j)}) - {\cal C}_{kj}
  - 4 \pi \sum_{1\leq i \leq N, ~ i \neq k}   {\cal C}_{ij} 
  ~\mbox{cap} (F^{(i)}) ~ G(\BO^{(k)}, \BO^{(i)}) =0, ~~ k,j = 1,\ldots, N,
 \eequ{T3_sys_a}
 and hence using \eq{T3_eq3}--\eq{T3_sys_a}, we arrive at
 \beq
 \CR(\Bx, \By) =  \CR^{(1)}(\Bx, \By)
 \eequ{T4_eq_w}
 for all $\Bx \in \prt ({\Bbb R}^3 \setminus F^{(k)})$ and $\By \in \GO_N$.  
 
 Let us consider the case when $\By \in \prt ({\Bbb R}^3 \setminus F^{(m)})$. 
Then
 $$
 T^{(j)}(\By) = 4 \pi ~\mbox{cap} (F^{(j)})~ G(\BO^{(j)}, \By) +O\Big(\fr{\Gve ~\mbox{cap} (F^{(j)}) }{| \By - \BO^{(j)} |^2}\Big), ~ j \neq m,
 $$
 and
 $$
  T^{(m)}(\By) = 1 -  4 \pi ~\mbox{cap} (F^{(m)})~ H(\BO^{(m)}, \By).
 $$
 The double sum in  \eq{T3_eq4} can be rearranged according to  \eq{T3_sys} 
\beqa
&&  \sum_{1\leq i \leq N, ~ i \neq k}    
   O \Big(  \fr{\Gve ~\mbox{cap} (F^{(i)})}{|  \BO^{(k)} - \BO^{(i)}|^2} \Big)   
   \sum_{j=1}^N
   {\cal C}_{ij} T^{(j)}(\By)  
  \nonumber \\
&&  =  \sum_{1\leq i \leq N, ~ i \neq k}    
   O \Big(  \fr{\Gve ~\mbox{cap} (F^{(i)})}{|  \BO^{(k)} - \BO^{(i)}|^2} \Big)    \Bigg\{ {\cal C}_{im} + O(|{\cal C}_{im}| \Gve)
   \nonumber \\
&&   +   4 \pi  \sum_{1 \leq j \leq N, j \neq m}  {\cal C}_{ij}  \Big( \mbox{ cap} (F^{(j)}) G(\BO^{(j)}, \By) + O(\fr{\Gve^2}{|\By - \BO^{(j)}|^2}) \Big) \Bigg\}
  \nonumber \\
&&  = \sum_{1\leq i \leq N, ~ i \neq k, i \neq m}    
   O \Big(  \fr{\Gve^2}{|  \BO^{(k)} - \BO^{(i)}|^2} \Big)  \Big\{ G(\BO^{(m)}, \BO^{(i)}) + O (\Gve d^{-1})  \Big \}
  \nonumber \\
  &&  +  \sum_{1\leq i \leq N, ~ i \neq k } ~ \sum_{1 \leq j \leq N, j \neq m}  {\cal C}_{ij} O(\fr{\Gve^4}{|\BO^{(m)} - \BO^{(j)}|^2  |  \BO^{(k)} - \BO^{(i)}|^2  })
   \l{T4_eq_fa} \\
&&   = O(\Gve^2 |\log d| d^{-3} + \Gve^3 d^{-4} + \Gve^4 d^{-7})
   , ~~ \mbox{for } \Bx \in \prt ( {\Bbb R}^3 \setminus F^{(k)}), ~\By \in \prt ( {\Bbb R}^3 \setminus F^{(m)}),
  \nonumber
  \end{eqnarray}
  where the estimate of the last double sum in   \eq{T4_eq_fa} is similar to     \eq{T3_eq1_a}. Combining  \eq{T3_eq4},  \eq{T4_eq_w} and   \eq{T4_eq_fa}, we deduce that
    $   \CR(\Bx, \By) = O(\Gve d^{-2}) $ for $\Bx \in \prt ( {\Bbb R}^3 \setminus F^{(k)}), ~\By \in \prt ( {\Bbb R}^3 \setminus F^{(m)}),$   $m,k = 1,\ldots,N.$
   Using the symmetry of $\CR(\Bx, \By)$ together with      \eq{T4_eq_g}  we also obtain that $\CR(\Bx, \By) = O(\Gve d^{-2}) $ for $\Bx \in \prt ( {\Bbb R}^3 \setminus F^{(k)}),    ~k = 1,\ldots,N,~\By \in 
   \prt \GO.$ Applying the maximum principle for harmonic functions we get
   \beq
   \CR(\Bx, \By) = O(\Gve d^{-2}) ~\mbox{for } \Bx \in \prt ( {\Bbb R}^3 \setminus F^{(k)}), k=1,\ldots,N,~\By \in \GO_N.
   \eequ{T4_eq_x}

Finally, formulae    \eq{T4_eq_h}    and \eq{T4_eq_x} imply that $\CR(\Bx, \By) = O(\Gve d^{-2})$ for $\Bx \in \prt \GO_N$ and $\By \in \GO_N$, 
and then applying the maximum principle for harmonic functions we complete the proof. $\Box$

 \vspace{.2in}
 
 \noindent {\bf Acknowledgement.} The support of the UK  Engineering and Physical Sciences Research Council via the grant EP/F005563/1 is gratefully acknowledged.  
 

 \end{document}